\documentclass[12pt,reqno,english, ]{article}
\usepackage[T1]{fontenc}
\usepackage[utf8]{inputenc}
\usepackage{geometry}
\usepackage{color}
\usepackage{babel}
\usepackage{amsmath}
\usepackage{amsthm}
\usepackage{amssymb}
\usepackage{graphicx}
\usepackage[unicode=true,pdfusetitle,
 bookmarks=true,bookmarksnumbered=false,bookmarksopen=false,
 breaklinks=false,pdfborder={0 0 1},backref=page,colorlinks=true]
 {hyperref}
\hypersetup{
 linkcolor=blue,citecolor=blue}

\makeatletter
\numberwithin{equation}{section}
\numberwithin{figure}{section}
\theoremstyle{plain}
\newtheorem{thm}{\protect\theoremname}[section]
\theoremstyle{plain}
\newtheorem{conjecture}[thm]{\protect\conjecturename}
\theoremstyle{definition}
\newtheorem*{defn*}{\protect\definitionname}
\theoremstyle{definition}
\newtheorem{defn}[thm]{\protect\definitionname}
\theoremstyle{plain}
\newtheorem{lem}[thm]{\protect\lemmaname}
\theoremstyle{plain}
\newtheorem{cor}[thm]{\protect\corollaryname}
\theoremstyle{remark}
\newtheorem{rem}[thm]{\protect\remarkname}
\theoremstyle{definition}
\newtheorem{example}[thm]{\protect\examplename}

\usepackage{ytableau}
\usepackage{appendix}
\usepackage[shortlabels]{enumitem}

\DeclareMathOperator{\Sym}{Sym}

\newlength{\tempindent} 
\newcommand{\lazyenum}{
\setlength{\tempindent}{\parindent} 
\begin{enumerate}[leftmargin=0cm,itemindent=0.7cm,labelwidth=\itemindent,labelsep=0cm,align=left,label=\arabic*)]
\setlength{\parskip}{\smallskipamount}
\setlength{\parindent}{\tempindent}
}

\makeatletter
\ifx\hyperlink\undefined
\newcommand{\customlabel}[2]{%
\protected@write \@auxout {}{\string \newlabel {#1}{{#2}{}}}}
\else
\newcommand{\customlabel}[2]{%
   \protected@write \@auxout {}{\string \newlabel {#1}{{#2}{\thepage}{#2}{#1}{}} }%
   \hypertarget{#1}%
}
\fi
\makeatother

\theoremstyle{plain}
\newtheorem*{thm:CLRp}{Theorem \ref{thm:CLR}'}
\newtheorem*{conj:Caputo's--shuffles-conjecturep}{Conjecture \ref{conj:Caputo's--shuffles-conjecture}'}
\newtheorem*{thm:triangles and the likep}{Theorem \ref{thm:triangles and the like}'}

\makeatother

\usepackage{listings}
\providecommand{\conjecturename}{Conjecture}
\providecommand{\corollaryname}{Corollary}
\providecommand{\definitionname}{Definition}
\providecommand{\examplename}{Example}
\providecommand{\lemmaname}{Lemma}
\providecommand{\remarkname}{Remark}
\providecommand{\theoremname}{Theorem}

\begin{document}
\global\long\def\ds{ \mathrm{DS}}%
 
\global\long\def\R{ \mathrm{\mathbf{R}}}%
 
\global\long\def\D{ \mathrm{\mathbf{D}}}%
\global\long\def\I{ \mathrm{\mathbf{I}}}%
 
\global\long\def\defi{ \stackrel{\mathrm{def}}{=}}%
 
\global\long\def\id{\mathrm{id}}%
 
\global\long\def\david{\mathrm{\star}}%
 
\global\long\def\tr{\mathrm{tr}}%
 
\global\long\def\ijlm{ij\ell m}%

\title{On the Aldous-Caputo Spectral Gap Conjecture for Hypergraphs}
\author{Gil Alon, Gady Kozma and Doron Puder}
\maketitle
\begin{abstract}
In their celebrated paper \cite{caputo2010proof}, Caputo, Liggett
and Richthammer proved Aldous' conjecture and showed that for an arbitrary
finite graph, the spectral gap of the interchange process is equal
to the spectral gap of the underlying random walk. A crucial ingredient
in the proof was the Octopus Inequality --- a certain inequality
of operators in the group ring $\mathbb{R}\left[\Sym_{n}\right]$
of the symmetric group. Here we generalize the Octopus Inequality
and apply it to generalizing the Caputo-Liggett-Richthammer Theorem
to certain hypergraphs, proving some cases of a conjecture of Caputo. 

\tableofcontents{}
\end{abstract}

\section{Introduction}

\subsection{The interchange process and random walk on graphs and hypergraphs}

Consider a finite undirected graph $G$ on $n$ vertices. In the interchange
process ($\mathrm{IP}$) on $G$, $n$ distinct balls are placed at
the vertices of $G$, one ball at every vertex. At every step of the
process, one picks an edge of $G$ uniformly at random, and interchanges
the two balls at the endpoints of the chosen edge. Every step is independent
of the other steps. This is a Markov chain with $n!$ states and accordingly
has a spectrum of $n!$ eigenvalues. All eigenvalues are real, contained
in $\left[-1,1\right]$, with a trivial eigenvalue at $1$ corresponding
to the uniform distribution on all possible states.

One may also consider a simpler process on $G$, which we refer to
as the random walk ($\mathrm{RW}$) on $G$. In this process there
is a single ball, which is located at one of the vertices. At each
step, one picks an edge of $G$ uniformly at random as before. If
the ball sits at an endpoint of this edge, it is moved to the other
endpoint, and otherwise it stays in its place. Note that unlike the
ordinary simple random walk, we do not necessarily pick an edge incident
to the ball. This Markov chain consists of $n$ possible states and
thus has a spectrum of $n$ eigenvalues. As before, they are all real
and contained in $\left[-1,1\right]$, and the uniform distribution
on the vertices gives a trivial eigenvalue at $1$. 

It is a simple observation (see Section \ref{sec:Preliminaries-representation theory of Sn})
that the spectrum of $\mathrm{RW}$ is contained in that of $\mathrm{IP}$,
so, in particular, $\lambda_{2}^{\mathrm{RW}}\left(G\right)\le\lambda_{2}^{\mathrm{IP}}\left(G\right)$,
where $\lambda_{i}$\marginpar{$\lambda_{i}$} is the $i^{\mathrm{th}}$
largest eigenvalue, counted with multiplicity. Around 1992, Aldous
conjectured that for any graph, the spectral gaps $\lambda_{1}-\lambda_{2}=1-\lambda_{2}$
of both processes are, in fact, identical. Some partial results were
proven along the following years (see \cite{caputo2010proof} for
a survey), but it was only established in full in 2010, by Caputo,
Liggett and Richthammer. Moreover, they established the conjecture
for any fixed probability distribution on the edges of the graph $G$.
\begin{thm}
\label{thm:CLR}\cite[Thm.~1.1]{caputo2010proof} Let $G$ be a finite
weighted graph with non-negative weights on the edges. Then 
\begin{equation}
\lambda_{2}^{\mathrm{RW}}\left(G\right)=\lambda_{2}^{\mathrm{IP}}\left(G\right),\label{eq:Aldous-result-IP-RW}
\end{equation}
where $\mathrm{RW}$ and $\mathrm{IP}$ are the random walk and interchange
process on $G$, respectively, with each edge chosen at every step
with probability proportional to its weight.
\end{thm}

The equality \eqref{eq:Aldous-result-IP-RW} is remarkable: although
the spectrum of $\mathrm{IP}$ is much larger -- of size $n!$ --
its highest non-trivial eigenvalue always comes from the small pool
of $n-1$ non-trivial eigenvalues of the $\mathrm{RW}$ process.

It is natural to wonder what possible generalizations of this intriguing
result may hold. One such generalization, sometimes called \emph{the
$\alpha$-shuffle conjecture}, was suggested by P.~Caputo (\cite[Conj.~3]{piras2010generalizations},\cite[P.~301]{cesi2016few},\cite[P.~78]{aldous2020life}).
It is a natural generalization to hypergraphs of the Aldous-Caputo-Liggett-Richthammer
Theorem. A finite weighted hypergraph consists of a finite set of
vertices together with a collection of distinct subsets of the edges,
called hyper-edges, and a non-negative weight for each hyper-edge.
So a weighted graph is the special case of a weighted hypergraph where
all hyper-edges are of size $2$. As before, at each step of the $\mathrm{IP}$
or the $\mathrm{RW}$ we may pick a random hyper-edge with probability
proportional to its weight. Once a certain hyper-edge is chosen, say
of size $d$, pick a uniformly random permutation of the $d$ vertices
contained in the hyper-edge, among all $d!$ possible permutations,
and move around accordingly the balls at these vertices. Here we allow
hyper-edges of cardinality $0$ and $1$, although they only add laziness
to the two processes.
\begin{conjecture}[Caputo's conjecture for hypergraphs]
\label{conj:Caputo's--shuffles-conjecture} Let $G$ be an arbitrary
finite weighted hypergraph with non-negative weights. Then 
\[
\lambda_{2}^{\mathrm{RW}}\left(G\right)=\lambda_{2}^{\mathrm{IP}}\left(G\right),
\]
where $\mathrm{RW}$ and $\mathrm{IP}$ are the random walk and interchange
process on $G$, respectively, defined by picking, at every step,
a random hyper-edge with probabilities proportional to the weights
and then a uniformly distributed permutation supported on the vertices
of this hyper-edge.
\end{conjecture}

As in the graph case, the equality $\lambda_{1}^{\mathrm{RW}}\left(G\right)=\lambda_{1}^{\mathrm{IP}}\left(G\right)$
and the inequality $\lambda_{2}^{\mathrm{RW}}\left(G\right)\le\lambda_{2}^{\mathrm{IP}}\left(G\right)$
are both easy observations -- see Lemma \ref{lem:aldous equivalent to that std rules}.
Note that Theorem \ref{thm:CLR} is the special case of Conjecture
\ref{conj:Caputo's--shuffles-conjecture} where the hypergraph is
restricted to be an ordinary graph. Indeed, although in the statement
of Conjecture \ref{conj:Caputo's--shuffles-conjecture} we add laziness
to the processes -- at every step, with probability $1/2$, we choose
the identity permutation -- this change shifts the entire spectrum
of both processes linearly by $\lambda\mapsto\frac{\lambda+1}{2}$,
so \eqref{eq:Aldous-result-IP-RW} in one setup is equivalent to the
same equality in the other setup.

The special case of Conjecture \ref{conj:Caputo's--shuffles-conjecture}
where the weight of an hyper-edge depends only on its size was recently
settled in \cite[Thm.~1.8]{bristiel2024entropy}\footnote{This special case can also be shown by simple induction using the
proof technique of \cite[Cor.~2.9]{parzanchevski2020aldous}.}. Here we generalize Theorem \ref{thm:CLR} to include more cases
from this conjecture:
\begin{thm}
\label{thm:triangles and the like} Let $G$ be a finite weighted
hypergraph with non-negative weights and $B$ a subset of its vertices.
Assume that every hyper-edge $E$ of $G$ satisfies $E\supseteq B$
and $\left|E\setminus B\right|\le2$. Then 
\[
\lambda_{2}^{\mathrm{RW}}\left(G\right)=\lambda_{2}^{\mathrm{IP}}\left(G\right),
\]
where the $\mathrm{IP}$ and $\mathrm{RW}$ processes are defined
as in Conjecture \ref{conj:Caputo's--shuffles-conjecture}.
\end{thm}

\begin{figure}
\begin{centering}
\includegraphics[viewport=50bp 0bp 400bp 250bp,scale=0.6]{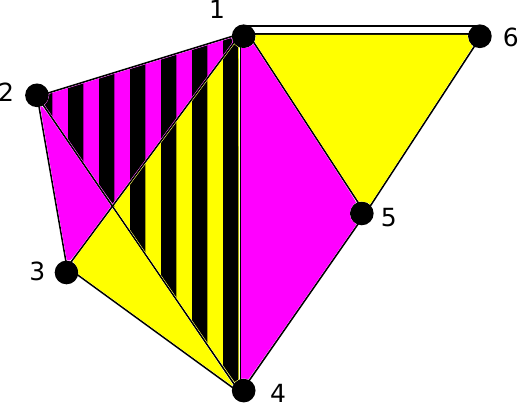}
\par\end{centering}
\caption{\label{fig:triangles}A hypergraph for which Theorem \ref{thm:triangles and the like}
applies. With vertices $\left[6\right]=\left\{ 1,2,3,4,5,6\right\} $,
this hypergraph consists of the hyper-edges $\left\{ 1,2,3\right\} $,
$\left\{ 1,2,4\right\} $, $\left\{ 1,3,4\right\} $, $\left\{ 1,4,5\right\} $,
$\left\{ 1,5,6\right\} $ and $\left\{ 1,6\right\} $. In the notation
of Theorem \ref{thm:triangles and the like}, here $B=\left\{ 1\right\} $.}
\end{figure}

When $B=\emptyset$, Theorem \ref{thm:triangles and the like} recovers
Theorem \ref{thm:CLR}. The setup in Theorem \ref{thm:triangles and the like}
can be thought of as a blow-up of the setup in Theorem \ref{thm:CLR}:
one takes the edges of the graph from Theorem \ref{thm:CLR}, adds
some singletons and possibly the empty set (both of which only add
harmless laziness in the original setup), and then adds the external
points of $B$ to each one of the sets. To illustrate, Theorem \ref{thm:triangles and the like}
applies to the hypergraph in Figure \ref{fig:triangles} (with arbitrary
non-negative weights on the hyper-edges).

We remark that limited simulations we conducted support Conjecture
\ref{conj:Caputo's--shuffles-conjecture} in its full generality.
We can also prove the conjecture in several additional cases not covered
by Theorem \ref{thm:triangles and the like}. These cases include
certain trees of hyper-edges or hypergraphs consisting of hyper-edges
of size at least $n-1$. We elaborate more in Section \ref{sec:Further-results-around-Caputo-conjecture}.

\subsection{A representation theoretic point of view}

A fruitful point of view on Aldous' spectral gap conjecture and the
results and conjectures above is that of representation theory. In
fact, this whole line of research began in the seminal paper of Diaconis-Shahshahani
\cite{diaconis1981generating}, which uses representation theory to
analyze the $\mathrm{IP}$ of the complete graph. In addition, Cesi
shows in \cite{cesi2016few} how representation theory clarifies the
proof of Aldous' conjecture from \cite{caputo2010proof}. In the current
paper this perspective is extremely useful and we stick to it.
\begin{defn*}
Denote the symmetric group on $n$ elements by $\mathrm{Sym}_{n}$
and denote $\left[n\right]=\left\{ 1,\ldots,n\right\} $. Consider
the group ring element
\begin{equation}
\Sigma=\sum_{1\le i<j\le n}c_{i,j}\left(i~j\right)\in\mathbb{R}\left[\Sym_{n}\right],\label{eq:aldous element}
\end{equation}
where $\left(i~j\right)$ is the transposition of $i$ and $j$ in
$\Sym_{n}$ and $c_{i,j}\ge0$ for all $1\le i<j\le n$. For every
finite-dimensional complex or real representation $\rho\colon\Sym_{n}\to\mathrm{GL}\left(V\right)$
of $\Sym_{n}$, the eigenvalues of $\rho\left(\Sigma\right)$ are
all real (see Footnote \ref{fn:real eigenvalues} below). We denote
by $\lambda_{i}\left(\Sigma,V\right)$\marginpar{$\lambda_{i}\left(\Sigma,V\right)$}
the $i$-th largest eigenvalue of $\rho\left(\Sigma\right)$, counted
with multiplicity. 

In these terms, the Caputo-Liggett-Richthammer theorem can be stated
as follows.
\end{defn*}
\begin{thm:CLRp}[Representation theoretic version of Theorem \ref{thm:CLR}]\customlabel{thm:CLR-rep-version}{\ref*{thm:CLR}'}
Let ${\bf R}=\mathbb{R}[\Sym_{n}]$ denote the regular representation
of $\Sym_{n}$ and let $\mathbf{D=\mathbb{R}}^{n}$ be the standard
representation of $\Sym_{n}$ via permutations of the coordinates.
For $\Sigma\in\mathbb{R}\left[\Sym_{n}\right]$ as in \eqref{eq:aldous element}
we have 
\begin{equation}
\lambda_{2}\left(\Sigma,\D\right)=\lambda_{2}\left(\Sigma,\R\right).\label{eq:Aldous result}
\end{equation}
\end{thm:CLRp}

In this setup, the eigenvalues of $\Sigma$ in both $\R$ and $\D$
are all contained in the interval $\left[-\sum c_{i,j},\text{\ensuremath{\sum c_{i,j}}}\right]$.
The equality $\lambda_{1}\left(\Sigma,\D\right)=\lambda_{1}\left(\Sigma,\R\right)=\lambda_{1}\left(\Sigma,\I\right)=\sum_{i,j}c_{i,j}$
(here $\I$ is the trivial representation), and the inequality $\lambda_{2}\left(\Sigma,\D\right)\le\lambda_{2}\left(\Sigma,\R\right)$
are straightforward: indeed, $\D$ is a sub-representation of $\R$,
and $\I$ is a subrepresentation of both --- see Lemma \ref{lem:aldous equivalent to that std rules}. 

It is natural to wonder which further elements $\Sigma\in\mathbb{R}\left[\Sym_{n}\right]$
satisfy the equality \eqref{eq:Aldous result} or similar equalities.
See the discussions in \cite[\S 5]{cesi2016few} and in \cite[\S 1 and \S 4]{parzanchevski2020aldous}
around this question. In this language, Caputo's $\alpha$-shuffles
conjecture becomes the following. For every $A\subseteq n$ define
a $J_{A}\in\mathbb{R}[\Sym_{n}]$ using\marginpar{$J_{A}$} 
\begin{equation}
{\displaystyle J_{A}\defi\sum_{\sigma\in\Sym_{n}\colon\mathrm{supp}\left(\sigma\right)\subseteq A}\sigma}.\label{eq:JA}
\end{equation}
For example, $J_{\left\{ 1,2,4\right\} }=\id+\left(1~2\right)+\left(1~4\right)+\left(2~4\right)+\left(1~2~4\right)+\left(1~4~2\right)$
and $J_{\left\{ 5\right\} }=J_{\emptyset}=\id$.

\begin{conj:Caputo's--shuffles-conjecturep}[Representation theoretic version of Conjecture \ref{conj:Caputo's--shuffles-conjecture}]\customlabel{conj:Caputo's--shuffles-conjecture-Rep-version}{\ref*{conj:Caputo's--shuffles-conjecture}'}Let
$\left\{ w_{A}:A\subseteq\left[n\right]\right\} $ be non-negative
real numbers, and consider the group ring element
\begin{equation}
U=\sum_{A\subseteq\left[n\right]}w_{A}J_{A}\in\mathbb{R}\left[\Sym_{n}\right].\label{eq:hypergraph element}
\end{equation}
Then $\lambda_{2}\left(U,{\bf D}\right)=\lambda_{2}\left(U,\mathbf{R}\right)$.

\end{conj:Caputo's--shuffles-conjecturep}

Theorem \ref{thm:triangles and the like} can be stated as follows.

\begin{thm:triangles and the likep}[Representation theoretic version of Theorem \ref{thm:triangles and the like}]
\customlabel{thm:triangles and the like - Rep version}{\ref*{thm:triangles and the like}'}Let $B\subseteq\left[n\right]$ and consider the group ring element
\[
U=\sum_{\substack{B\subseteq A\subseteq\left[n\right]\colon\\
\left|A\setminus B\right|\le2
}
}w_{A}J_{A}\in\mathbb{R}\left[\Sym_{n}\right],
\]
where $w_{A}\ge0$ are arbitrary non-negative numbers. Then $\lambda_{2}\left(U,{\bf D}\right)=\lambda_{2}\left(U,\mathbf{R}\right)$.\end{thm:triangles and the likep}

\subsection{Generalizations of the Octopus Inequality}

A key ingredient of the proof of Theorem \ref{thm:CLR-rep-version}
in \cite{caputo2010proof}, which we sketch in Section \ref{subsec:A-sketch-of-CLR-proof},
is a certain inequality of operators in the group ring of $\Sym_{n}$,
called the Octopus Inequality. The Octopus Inequality is also of independent
interest, see \cite{chen2017moving} and \cite{alon2020comparing}.
Inspired by the proof in \cite{caputo2010proof}, a main ingredient
of our proof of Theorem \ref{thm:triangles and the like - Rep version}
is a generalization of the Octopus Inequality. In fact, we also include
further generalizations of the Octopus Inequality which are not used
in the proof of Theorem \ref{thm:triangles and the like - Rep version}.
We find this type of general inequalities of operators in the group
algebra $\mathbb{R}\left[\Sym_{n}\right]$ remarkable and interesting
for their own right. 
\begin{defn}
\label{def:inequality of group algebra elements}For two elements
$\Sigma_{1},\Sigma_{2}\in\mathbb{R}\left[\Sym_{n}\right]$ we write
$\Sigma_{1}\ge\Sigma_{2}$ if $\rho\left(\Sigma_{1}-\Sigma_{2}\right)$
is a positive semi-definite operator for every (finite dimensional)
representation $\rho$ of $\Sym_{n}$.\footnote{\label{fn:positive-semidefinite}For the representation $\rho\colon\Sym_{n}\to\mathrm{GL}\left(V\right)$,
the operator $M\defi\rho\left(\Sigma_{1}-\Sigma_{2}\right)$ is called
positive semi-definite if for every $v\in V$ we have $\left\langle Mv,v\right\rangle \ge0$
in the standard inner product on $V$. This is equivalent to $M$
being Hermitian and having non-negative eigenvalues.} 
\end{defn}

\begin{defn}
For any subset $\emptyset\ne A\subseteq\left[n\right]$ denote\marginpar{$\alpha_{A}$}
\begin{equation}
\alpha_{A}\defi\frac{1}{\left(\left|A\right|-1\right)!}\sum_{\pi\in\mathrm{Sym}\left(A\right)}\left(\id-\pi\right)\in\mathbb{R}\left[\Sym_{n}\right],\label{eq:def of alpha}
\end{equation}
where $\id\in\Sym_{n}$ denotes the identity permutation. 
\end{defn}

Thus $\alpha_{A}=\left|A\right|\cdot\id-\frac{J_{A}}{\left(\left|A\right|-1\right)!}$.
For example, $\alpha_{\left\{ i,j\right\} }=\id-\left(i~j\right)$.
For completeness, we also define $\alpha_{\emptyset}\defi0\in\mathbb{R}\left[\Sym_{n}\right]$
for $A$ the empty set. Note that $\alpha_{A}=0$ whenever $\left|A\right|\le1$.
As explained in Lemma \ref{lem:alpha is non-negative}, $\alpha_{A}\ge0$
for all $A\subseteq\left[n\right]$.

The Octopus Inequality \cite[Thm.~2.3]{caputo2010proof} (and see
also \cite[Eq.~(4.6)]{cesi2016few}), can be formulated as follows:
\begin{thm}[Octopus Inequality]
\label{thm:octopus} Let $c_{2},\ldots,c_{n}\ge0$ be non-negative
numbers. Then we have the following inequality of operators in $\mathbb{R}\left[\Sym_{n}\right]$:
\begin{equation}
\left(\sum_{i=2}^{n}c_{i}\right)\left(\sum_{j=2}^{n}c_{j}\alpha_{\left\{ 1,j\right\} }\right)\ge\sum_{2\le i<j\le n}c_{i}c_{j}\alpha_{\left\{ i,j\right\} }.\label{eq:octopus}
\end{equation}
\end{thm}

One generalization of the Octopus Inequality is the following:
\begin{thm}[Generalized octopus --- disjoint sets]
\label{thm:squid - disjoint sets} Let $A_{1},\ldots,A_{t}\subseteq\left\{ 2,\ldots,n\right\} $
be pairwise \textbf{\emph{disjoint}} subsets and $c_{1},\ldots,c_{t}\ge0$
be non-negative numbers. Then we have the following inequality of
operators in $\mathbb{R}\left[\Sym_{n}\right]$:
\begin{equation}
\left(\sum_{i=1}^{t}c_{i}\left|A_{i}\right|\right)\left(\sum_{j=1}^{t}c_{j}\left(\alpha_{A_{j}\cup\left\{ 1\right\} }-\alpha_{A_{j}}\right)\right)\ge\sum_{i=1}^{t}c_{i}^{~2}\alpha_{A_{i}}+\sum_{1\le i<j\le t}c_{i}c_{j}\left(\alpha_{A_{i}\cup A_{j}}-\alpha_{A_{i}}-\alpha_{A_{j}}\right).\label{eq:squid - disjoint sets}
\end{equation}
\end{thm}

When all the $A_{i}$'s are singletons, Theorem \ref{thm:squid - disjoint sets}
recovers the Octopus Inequality (recall that $\alpha_{A}=0$ for $A$
a singleton). Yet another generalization of Theorem \ref{thm:octopus}
in a different direction is the following result, which is the key
ingredient in the proof of Theorem \ref{thm:triangles and the like - Rep version}.
\begin{thm}[Generalized octopus --- sets with large intersection]
\label{thm:squid - sets of size k containing a common (k-1)-subset}
Let $A_{0}\subseteq\left\{ 2,\ldots,n\right\} $, let $A_{1},\ldots,A_{t}\subseteq\left\{ 2,\ldots,n\right\} $
be distinct sets satisfying $A_{0}\subseteq A_{i}\subseteq\left\{ 2,\ldots,n\right\} $
and $\left|A_{i}\setminus A_{0}\right|=1$, and let $c_{0},c_{1},\ldots,c_{t}\ge0$
be non-negative numbers. Then we have the following inequality of
operators in $\mathbb{R}\left[\Sym_{n}\right]$:
\begin{equation}
\left(\sum_{i=0}^{t}c_{i}\left|A_{i}\right|\right)\left(\sum_{j=0}^{t}c_{j}\left(\alpha_{A_{j}\cup\left\{ 1\right\} }-\alpha_{A_{j}}\right)\right)\ge\sum_{i=0}^{t}c_{i}^{~2}\alpha_{A_{i}}+\sum_{0\le i<j\le t}c_{i}c_{j}\left(\alpha_{A_{i}\cup A_{j}}+\alpha_{A_{0}}\right).\label{eq:squid - sets with large intersection}
\end{equation}
\end{thm}

Theorem \ref{thm:squid - sets of size k containing a common (k-1)-subset},
too, restricts to Theorem \ref{thm:octopus}: this time, when $A_{0}=\emptyset$.
A third inequality we prove in this paper is the content of the following
theorem. This one does not contain the original Octopus Inequality
as a special case.
\begin{thm}[Generalized Octopus --- sets of co-size one]
\label{thm:squid - sets of cosize 1 in A0} Let $A_{0}\subseteq\left\{ 2,\ldots,n\right\} $,
let $A_{1},\ldots,A_{t}$ be distinct sets satisfying $A_{i}\subseteq A_{0}$
and $\left|A_{0}\setminus A_{i}\right|=1$ for $i=1,\ldots,t$, and
let $c_{0},c_{1},\ldots,c_{t}\ge0$ be non-negative numbers. Then
we have the following inequality of operators in $\mathbb{R}\left[\Sym_{n}\right]$:
\begin{equation}
\left(\sum_{i=0}^{t}c_{i}\left|A_{i}\right|\right)\left(\sum_{j=0}^{t}c_{j}\left(\alpha_{A_{j}\cup\left\{ 1\right\} }-\alpha_{A_{j}}\right)\right)\ge\sum_{i=0}^{t}c_{i}^{~2}\alpha_{A_{i}}+\sum_{0\le i<j\le t}c_{i}c_{j}\left(\alpha_{A_{0}}+\alpha_{A_{i}\cap A_{j}}\right).\label{eq:squid - sets of cosize-1 in A0}
\end{equation}
\end{thm}

In Section \ref{subsec:A-general-inequality} we explain how the three
Theorems \ref{thm:squid - disjoint sets}, \ref{thm:squid - sets of size k containing a common (k-1)-subset}
and \ref{thm:squid - sets of cosize 1 in A0} all originate from the
same inequality scheme. 

\subsection*{Paper organization }

We begin in Section \ref{sec:Preliminaries-representation theory of Sn}
with some preliminaries on the representation theory of the symmetric
group, as well as some lemmas that will be used in the following sections.
Section \ref{sec:The-proof-of-Aldous} sketches the proof of Aldous'
conjecture by Caputo-Liggett-Richthammer, explains how one might generalize
it to hypergraphs, and reduces Theorem \ref{thm:triangles and the like - Rep version}
to Theorem \ref{thm:squid - sets of size k containing a common (k-1)-subset}.
Sections \ref{sec:Proof-of- disjoint sets}, \ref{sec:sets with large intersection}
and \ref{sec:n-1 sets} prove Theorems \ref{thm:squid - disjoint sets},
\ref{thm:squid - sets of size k containing a common (k-1)-subset}
and \ref{thm:squid - sets of cosize 1 in A0}, respectively, the proof
of Theorem \ref{thm:squid - sets of size k containing a common (k-1)-subset}
being by far the most involved. In fact, in the latter proof we also
rely partially on computer computations. In Section \ref{sec:Further-results-around-Caputo-conjecture}
we mention a few further special cases of Caputo's conjecture on hypergraphs
which hold true, and in Section \ref{sec:Counterexamples} mention
a few examples illustrating the limitations of our approach as detailed
in Section \ref{sec:The-proof-of-Aldous}. Finally, the appendix contains
the sage-code we used in some arguments in the proof of Theorem \ref{thm:squid - sets of size k containing a common (k-1)-subset},
as detailed in Section \ref{sec:sets with large intersection}. 

\subsection*{Notation}

Throughout the paper we repeatedly use the notation of $\alpha_{A}$
as defined in \eqref{eq:def of alpha} and of $J_{A}$ as defined
in \eqref{eq:JA}. If $E$ is a (finite) set and $d\in\mathbb{Z}_{\ge0}$
a non-negative integer, we denote by $\binom{E}{d}$ the set of subsets
of $E$ of size $d$. 

\subsection*{Acknowledgments}

We thank the anonymous referee for suggesting some improvements to
a previous version. G.K.~and D.P.~were supported by the Israel Science
Foundation, ISF grants 607/21 and 1071/16, respectively. D.P.~was
also supported by the European Research Council (ERC) under the European
Union’s Horizon 2020 research and innovation programme (grant agreement
No 850956). 

\section{Representation theory preliminaries \label{sec:Preliminaries-representation theory of Sn}}

Let us start by a quick overview of the representation theory of the
symmetric group. As in the introduction, we set ${\bf R}=\mathbb{R}[\Sym_{n}]$,
$\mathbf{D=\mathbb{R}}^{n}$ and $\I=\mathbb{R}$, where all three
are considered as representations of $\Sym_{n}$. We shall often use
the ring embedding $\mathbb{R}\hookrightarrow\mathbf{R}$ via $c\mapsto c\cdot\id$,
in an implicit manner.

\begin{figure}
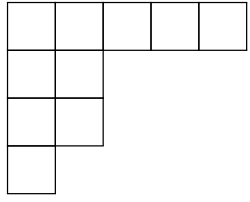

\begin{centering}
\ydiagram{5,2,2,1}
\par\end{centering}
\caption{\label{fig:Young-diagram}The Young diagram corresponding to the partition
$\mu=\left(5,2,2,1\right)\vdash10$}
\end{figure}

Recall that a partition $\mu=\left(\mu_{1},\dotsc,\mu_{r}\right)$
of a positive integer $n$ is a sequence of integers $\mu_{1}\geq\mu_{2}\geq\dotsb\geq\mu_{r}>0$
such that $\sum\mu_{i}=n$. This occurrence is denoted by \marginpar{$\mu\vdash n$}$\mu\vdash n$.
We also set $|\mu|=n$. A \emph{Young Diagram} is a finite array of
squares, which is made out of left aligned rows, such that the lengths
of the rows (when viewed from top to bottom) are non-increasing. Therefore,
there is a one-to-one correspondence between partitions of $n$ and
Young diagrams of size $n$: a partition $\mu$ corresponds to a Young
diagram whose $i^{\mathrm{th}}$ row has $\mu_{i}$ squares. See Figure
\ref{fig:Young-diagram}. In the sequel, we shall not distinguish
between a partition $\mu$ and its corresponding Young diagram.

As a finite group, the symmetric group $\Sym_{n}$ has only a finite
number of isomorphism classes of representations over $\mathbb{C}$
(this number is equal to the number of conjugacy classes of the group),
and in the case of $\Sym_{n}$ they can all be realized over\footnote{This and all other standard facts about the representation theory
of $\Sym_{n}$ can be found in any standard book on the subject, such
as \cite{fulton2013representation}.} $\mathbb{R}$. These isomorphism classes are classified by the Young
diagrams of size $n$. For any $\mu\vdash n$, we let $V_{\mu}$\marginpar{$V_{\mu}$}
denote a representation of $\Sym_{n}$ whose isomorphism class corresponds
to $\mu$. The following facts are well known:
\begin{itemize}
\item $\I\cong V_{\left(n\right)}$ is the trivial (one-dimensional) representation
of $\Sym_{n}$,
\item $\mathbf{D}\cong V_{\left(n\right)}\oplus V_{\left(n-1,1\right)}$
is the standard representation of $\Sym_{n}$, and
\item ${\bf R}\cong\bigoplus_{\mu\vdash n}\dim V_{\mu}\cdot V_{\mu}$ is
the regular representation of $\Sym_{n}$.
\end{itemize}
The last isomorphism means that the decomposition of the regular representation
$\R$ into irreducible representations contains all the irreducible
ones, each with multiplicity equal to its dimension. In the sequel,
we denote these three representations $\I_{n}$, $\D_{n}$ and $\R_{n}$
when the parameter $n$ is not clear from the context. 

We now expand the notation from Definition \ref{def:inequality of group algebra elements}.
For any representation $V$ of $\Sym_{n}$ and any group ring element
$\Sigma\in\mathbb{R}\left[\Sym_{n}\right]$, we shall denote by $\Sigma|_{V}$\marginpar{$\Sigma|_{V}$}
the action of $\Sigma$ on $V$, which is an element of $\mathrm{End}\left(V\right)$.
If $\Sigma$ is symmetric (i.e., $\Sigma=\sum_{\pi\in\Sym_{n}}a_{\pi}\pi$
with $a_{\pi}=a_{\pi^{-1}}$ for all $\pi$), then the eigenvalues
of $\Sigma|_{V}$ are real\footnote{\label{fn:real eigenvalues}Indeed, every $d$-dimensional complex
representation $\rho\colon G\to\mathrm{GL}\left(V\right)$ of a finite
group $G$ admits a basis in which $\rho\left(g\right)\in U\left(d\right)$
is unitary for all $g\in G$. If $\Sigma\in\mathbb{C}\left[G\right]$
is symmetric, then $\rho\left(\Sigma\right)$ is Hermitian and so
has real eigenvalues.}. For any $\Sigma_{1},\Sigma_{2}\in\mathbb{R}\left[\Sym_{n}\right]$,
we shall denote by $\Sigma_{1}=_{V}\Sigma_{2}$\marginpar{${\scriptstyle \Sigma_{1}=_{V}\Sigma_{2}}$}
the fact that $\left(\Sigma_{1}-\Sigma_{2}\right)|_{V}=0$, and by
$\Sigma_{1}\geq_{V}\Sigma_{2}$\marginpar{${\scriptstyle \Sigma_{1}\protect\geq_{V}\Sigma_{2}}$}
(respectively, $\Sigma_{1}>_{V}\Sigma_{2}$) the fact that $\left(\Sigma_{1}-\Sigma_{2}\right)|_{V}$
is a positive semi-definite (respectively, positive definite) operator
(see Footnote \ref{fn:positive-semidefinite}). If $\Sigma_{1}-\Sigma_{2}$
is symmetric, then this is equivalent to $\Sigma_{1}-\Sigma_{2}$
having only non-negative (respectively, positive) eigenvalues. Whenever
$\Sigma_{1}\geq_{\mathbf{R}}\Sigma_{2}$ we shall simply write $\Sigma_{1}\geq\Sigma_{2}$\marginpar{${\scriptstyle \Sigma_{1}\protect\geq\Sigma_{2}}$}.
Note that this is equivalent to $\Sigma_{1}\geq_{V_{\mu}}\Sigma_{2}$
holding for every $\mu\vdash n$. 

The following basic lemma will be very useful.
\begin{lem}
\label{lem:aldous equivalent to that std rules}Let $\Sigma=\sum_{\pi\in\Sym_{n}}a_{\pi}\pi\in\mathbb{R}\left[\Sym_{n}\right]$
be a symmetric element with $a_{\pi}\ge0$ for all $\pi$ (in particular,
as stated above, it has a real spectrum in every representation).
Then $\lambda_{1}(\Sigma,\D)=\lambda_{1}(\Sigma,\R)=\lambda_{1}(\Sigma,\I)=\sum_{\pi}a_{\pi}$
and $\lambda_{2}(\Sigma,\D)\le\lambda_{2}(\Sigma,\R)$. Moreover,
the statement $\lambda_{2}\left(\Sigma,\D\right)=\lambda_{2}\left(\Sigma,\R\right)$
(as in Conjecture \ref{conj:Caputo's--shuffles-conjecture-Rep-version})
is equivalent to that for every $\mu\vdash n$, $\mu\ne\left(n\right),(n-1,1)$,
\[
\lambda_{1}\left(\Sigma,V_{\mu}\right)\le\lambda_{1}\left(\Sigma,V_{\left(n-1,1\right)}\right).
\]
\end{lem}

\begin{proof}
For every $\pi\in\Sym_{n}$ and every finite dimensional representation
$\rho$, the eigenvalues of $\rho(\pi)$ are roots of unity so, in
particular, $\left\Vert \rho(\pi)+\rho(\pi^{-1})\right\Vert \le2$.
Hence the eigenvalue of the trivial representation of $\Sigma$, which
is $\sum_{\pi\in\Sym_{n}}a_{\pi}$, is at least as large as any eigenvalue
of any other representation. As $\I$ is contained in both $\D$ and
$\R$, we obtain $\lambda_{1}(\Sigma,\D)=\lambda_{1}(\Sigma,\R)=\lambda_{1}(\Sigma,\I)=\sum_{\pi}a_{\pi}$.
The other statements of the lemma follow from the decompositions,
mentioned above, of $\D$ and $\R$ to irreducible representations.
\end{proof}
Every finite-dimensional representation $V$ of $\Sym_{n}$ admits
a $\Sym_{n}$-invariant inner product. If $\Sigma\in\mathbb{R}\left[\Sym_{n}\right]$
is symmetric, then we have a spectral decomposition $V=\bigoplus V_{\lambda}$
where $V_{\lambda}=\ker\left(\Sigma-\lambda\right)$ is the eigenspace
of the eigenvalue $\lambda$. For every $\lambda_{1}\ne\lambda_{2}$,
$V_{\lambda_{1}}$ and $V_{\lambda_{2}}$ are orthogonal with respect
to the inner product on $V$. Consequently, 
\begin{equation}
V_{\lambda_{1}}^{\bot}=\bigoplus_{\lambda\ne\lambda_{1}}V_{\lambda},\label{eq:orthogonal complement}
\end{equation}
regardless of which inner product we use. 

The following lemmas are needed in the next sections.
\begin{lem}
\label{lem:positive_on_Sk}Let $A\subseteq\left[n\right]$ and let
$\Sigma$ be in the group ring $\mathbb{R}[\Sym_{A}]$. Then $\Sigma|_{\mathbb{R}[\Sym_{A}]}\geq0$
($\Sigma|_{\mathbb{R}[\Sym_{A}]}=0$) if and only if $\Sigma|_{\mathbf{R}}\geq0$
($\Sigma|_{\mathbf{R}}=0$, respectively), where $\Sigma$ is viewed
as an element of $\mathbf{R}=\mathbb{R}\left[\Sym_{n}\right]$ through
the natural embedding $\Sym_{A}\subseteq\Sym_{n}$.
\end{lem}

\begin{proof}
Let $g_{1},\dotsc,g_{k}\in\Sym_{n}$ be representatives of the right
cosets of $\Sym_{A}$ in $\Sym_{n}$. Each $g_{i}$ generates in $\mathbf{R}$
a $\Sym_{A}$-submodule which is isomorphic, as a $\Sym_{A}$-submodule,
to $\mathbb{R}[\Sym_{A}]$. As $\mathbf{R}$ decomposes as a direct
sum of these submodules, the lemma follows.
\end{proof}
Let $A\subseteq\left[n\right]$ be any subset. Recall the notations
$J_{A}=\sum_{\pi\in\Sym_{A}}\pi\in\mathbb{R}\left[\Sym_{n}\right]$
from \eqref{eq:JA}, and $\alpha_{A}=\frac{1}{\left(\left|A\right|-1\right)!}\sum_{\pi\in\Sym_{A}}\left(\id-\pi\right)\in\mathbb{R}\left[\Sym_{n}\right]$
from \eqref{eq:def of alpha} (with $\alpha_{\emptyset}\defi0$). 
\begin{lem}
\label{lem:alpha is non-negative}For any $A\subseteq[n]$, $\alpha_{A}\ge0$.
\end{lem}

\begin{proof}
By Lemma \ref{lem:aldous equivalent to that std rules}, $\lambda_{1}(J_{A},\R)=\left|A\right|!$.
The lemma follows as $\alpha_{A}=|A|\cdot\id-\frac{J_{A}}{\left(\left|A\right|-1\right)!}$.
\end{proof}
\begin{lem}
\label{lem:alpha_A in S_A}Let $\mu\vdash n$. Then 
\[
\alpha_{\left[n\right]}\Big|_{V_{\mu}}\equiv\begin{cases}
0 & \mathrm{if}~\mu=\left(n\right),\\
n & \mathrm{otherwise.}
\end{cases}
\]
\end{lem}

\begin{proof}
Recall that $J_{\left[n\right]}=\sum_{\pi\in\Sym_{n}}\pi$. Since
$J_{[n]}$ commutes with every element of $\Sym_{n}$, by Schur's
lemma, $J_{[n]}$ acts as a scalar on any irreducible representation
of $\Sym_{n}$. Hence $J_{[n]}|_{V_{\mu}}=c$ for some $c\in\mathbb{R}$.
Since $\pi J_{[n]}=J_{[n]}$ for all $\pi\in\Sym_{n}$, if $V_{\mu}$
is non-trivial then $c=0$. If $\mu=(n)$, namely, if $V_{\mu}$ is
the trivial representation, then $J_{[n]}|_{V_{\mu}}=n!$. The result
now follows from $\alpha_{\left[n\right]}=\frac{1}{(n-1)!}\left(n!-J_{[n]}\right)$.
\end{proof}
\begin{cor}
\label{cor:alpha is constant when many squares outside first row}Let
$A\subseteq\left[n\right]$, and let $\mu\vdash n$ be a Young diagram
such that $\mu_{1}<|A|.$ Then $\alpha_{A}|_{V_{\mu}}\equiv\left|A\right|$.
\end{cor}

\begin{proof}
By the branching rule \cite[Exer.~4.53]{fulton2013representation},
the action of $\alpha_{A}$ on $V_{\mu}$ can be decomposed into the
action on irreducible representations of $\Sym_{\left|A\right|}$
with Young diagrams obtained from $\mu$ by removing $n-\left|A\right|$
squares (we apply the branching rule $n-|A|$ times). Since $\sum_{j\geq2}\mu_{j}>n-|A|$,
all these diagrams have more than one row, so the corresponding representations
of $\Sym_{A}$ are non-trivial. We are done by Lemma \ref{lem:alpha_A in S_A}.
\end{proof}
\begin{lem}
\label{lem:alpha_on_standard}For any $A\subseteq\left[n\right]$,
we have 
\begin{equation}
\alpha_{A}=_{\mathbf{D}}\sum_{\left\{ i,j\right\} \in\binom{A}{2}}\alpha_{\{i,j\}},\label{eq:alpha in D}
\end{equation}
and 
\begin{equation}
\alpha_{A}\le\sum_{\left\{ i,j\right\} \in\binom{A}{2}}\alpha_{\{i,j\}}.\label{eq:alpha_A le sum of alpha of pairs}
\end{equation}
\end{lem}

\begin{proof}
For any $\pi\in\Sym_{n}$, 
\[
\sum_{i\in\left\{ 1,\ldots,n\right\} }\left(i~~\pi(i)\right)=_{\mathbf{D}}\pi+\pi^{-1}+\left(n-2\right)
\]
(where if $\pi\left(i\right)=i$, the notation $\left(i~~\pi\left(i\right)\right)=\left(i~~i\right)$
represents the identity permutation). Indeed, the above equality can
be verified by comparing the action of the group ring elements on
both sides on the elements of the standard basis of $\mathbf{D}$.
Hence, 
\[
\sum_{i\in\left\{ 1,\ldots,n\right\} }\left(1-\left(i~~\pi(i)\right)\right)=_{\mathbf{D}}2-\pi-\pi^{-1},
\]
so
\[
\sum_{i\in\left\{ 1,\ldots,n\right\} }\alpha_{\left\{ i,\pi(i)\right\} }=_{\mathbf{D}}\left(1-\pi\right)+\left(1-\pi^{-1}\right).
\]
Summing over $\pi\in\Sym_{A}$, we get 

\[
\left(|A|-1\right)!\sum_{i,j\in A}\alpha_{\left\{ i,j\right\} }=_{\D}2\left(|A|-1\right)!\alpha_{A},
\]
which yields \eqref{eq:alpha in D}. 

For the inequality \eqref{eq:alpha_A le sum of alpha of pairs}, Lemma
\ref{lem:positive_on_Sk} shows that it is enough to show that $\alpha_{A}\le\sum_{\left\{ i,j\right\} \in\binom{A}{2}}\alpha_{\{i,j\}}$
in $\mathbb{R}\left[\Sym_{A}\right]$. Denote $k=\left|A\right|$.
Both sides are equal in the standard representation $\D_{k}=V_{\left(k\right)}\oplus V_{\left(k-1,1\right)}$
by the first part of the current lemma. We claim there is an inequality
in $V_{\mu}$ for every $\mu\vdash k$ with $\mu\ne\left(k\right),\left(k-1,1\right)$.
Indeed, the left hand side ($\alpha_{A}$) is equivalent to $\left|A\right|$
in $V_{\mu}$ for every $\left(k\right)\ne\mu\vdash k$ by Lemma \ref{lem:alpha_A in S_A}.
The right hand side $\sum_{\left\{ i,j\right\} \in\binom{A}{2}}\alpha_{\{i,j\}}$
is a conjugation-invariant element in $\mathbb{R}\left[\Sym_{A}\right]$
and thus by Schur's lemma acts as a scalar in every irreducible representation
$V_{\mu}$. Denoting by $\chi_{\mu}$ the character corresponding
to $V_{\mu}$, this scalar is
\begin{equation}
\binom{k}{2}\left(1-\frac{\chi_{\mu}\left(\left(1~2\right)\right)}{\chi_{\mu}\left(\id\right)}\right),\label{eq:DS-scalar}
\end{equation}
and \cite[Lemma 10]{diaconis1981generating} states that among all
$\mu\ne\left(k\right)$, the scalar \eqref{eq:DS-scalar} is smallest
for $\mu=\left(k-1,1\right)$. We conclude that because there is equality
$\alpha_{A}=_{V_{(k-1,1)}}\sum_{\left\{ i,j\right\} \in\binom{A}{2}}\alpha_{\{i,j\}}$,
there must be an inequality 
\[
\alpha_{A}\le_{V_{\mu}}\sum_{\{i,j\}\in\binom{A}{2}}\alpha_{\left\{ i,j\right\} }
\]
for every $\mu\vdash k$ with $\mu\ne\left(k\right),\left(k-1,1\right)$.
\end{proof}
We denote by $\tr$\marginpar{$\protect\tr$} the trace map.
\begin{lem}
\label{lem:The-trace-of-alpha_A in std} For any $A\subseteq\left[n\right]$,
we have $\tr\left(\alpha_{A}|_{\mathbf{D}}\right)=|A|\left(|A|-1\right)$.
\end{lem}

\begin{proof}
For $1\leq i<j\leq n$ we have $\alpha_{\{i,j\}}=1-(ij)$. Consider
the standard basis of $\D$ corresponding to the elements $1,\ldots,n$,
so that if $v=\left(v_{1},\ldots,v_{n}\right)\in\D$, then $\pi.v=\left(v_{\pi^{-1}\left(1\right)},\ldots,v_{\pi^{-1}\left(n\right)}\right)$
for all $\pi\in\Sym_{n}$. In this basis, the operator $\alpha_{\left\{ i,j\right\} }|_{\D}$
corresponds to the $n\times n$ matrix with all zeros except for the
$\left\{ i,j\right\} $-minor which is 
\[
\left(\begin{array}{cc}
1 & -1\\
-1 & 1
\end{array}\right),
\]
and we get $\tr\left(\alpha_{\{i,j\}}|_{\mathbf{D}}\right)=2$. The
result now follows from Lemma \ref{lem:alpha_on_standard}.
\end{proof}
\begin{lem}
\label{lem:if A subset of B}If $A\subseteq B\subseteq\left[n\right]$
then
\begin{enumerate}
\item \label{enu:product alpha_a alpha_b}$\alpha_{A}\alpha_{B}=\alpha_{B}\alpha_{A}=|B|\alpha_{A}$,
\item \label{enu:kernels}$\ker\alpha_{B}\subseteq\ker\alpha_{A}$ (here
$\alpha_{A}$ and $\alpha_{B}$ are viewed as linear operators on
${\bf R}=\mathbb{R}[\Sym_{n}]$),
\item \label{enu:eigenvalues of alpha_a}Each eigenvalue of $\alpha_{A}$
on $\mathbb{R}\left[\Sym_{n}\right]$ is either $0$ or $\left|A\right|$,
and
\item \label{enu:norm-of-alpha_a}$\left\Vert \alpha_{A}-\frac{|A|}{2}\right\Vert =\frac{|A|}{2}$
where $\left\Vert \cdot\right\Vert $ is the operator norm.
\end{enumerate}
\end{lem}

\begin{proof}
$\,$

\lazyenum \item

The claim is trivial when $A=\emptyset$ (as $\alpha_{\emptyset}=0$),
so assume $A\ne\emptyset$. Recall that $\alpha_{A}=|A|-\frac{J_{A}}{\left(\left|A\right|-1\right)!}$
and that $\Sym_{A}$ and $\Sym_{B}$ are naturally embedded in $\Sym_{n}$.
For each $\pi\in\Sym_{B}$ we have $\pi J_{B}=J_{B}\pi=J_{B}$. Since
$\Sym_{A}\subseteq\Sym_{B}$, we have 
\[
J_{B}J_{A}=J_{A}J_{B}=|A|!J_{B}.
\]
Hence, 
\begin{align*}
\alpha_{B}\alpha_{A} & =\alpha_{A}\alpha_{B}=\alpha_{A}\left(|B|-\frac{J_{B}}{\left(\left|B\right|-1\right)!}\right)\\
 & =|B|\alpha_{A}-\frac{\alpha_{A}J_{B}}{\left(\left|B\right|-1\right)!}=|B|\alpha_{A}-\frac{\left(|A|-\frac{J_{A}}{\left(\left|A\right|-1\right)!}\right)J_{B}}{\left(\left|B\right|-1\right)!}.
\end{align*}
Since $\left(|A|-\frac{J_{A}}{\left(\left|A\right|-1\right)!}\right)J_{B}=|A|\cdot J_{B}-\frac{|A|!J_{B}}{\left(\left|A\right|-1\right)!}=0$,
we conclude the desired equality $\alpha_{A}\alpha_{B}=\alpha_{B}\alpha_{A}=|B|\alpha_{A}$.

\item Let $v\in\ker\alpha_{B}$. By Item \ref{enu:product alpha_a alpha_b},
\[
\alpha_{A}v=\frac{1}{|B|}\alpha_{A}\left(\alpha_{B}v\right)=0.
\]

\item Note that $\alpha_{A}$ is a symmetric element of $\mathbb{R}\left[\Sym_{n}\right]$
and so is diagonalizable (with real eigenvalues). Setting $A=B$ in
Item \ref{enu:product alpha_a alpha_b}, we get $\alpha_{A}^{2}=|A|\alpha_{A}$,
from which the claim follows.

\item This follows immediately from Item \ref{enu:eigenvalues of alpha_a}.\qedhere\end{enumerate}
\end{proof}
Next we note the following observation regarding the ``Octopus''
element -- the difference between the two sides of \eqref{eq:octopus}.
This observation is basically contained in \cite[\S2]{caputo2010proof}.
\begin{lem}
\label{lem:Octopus is rank 1 in standard}Let 
\[
\Delta\defi\left(\sum_{i=2}^{n}c_{i}\right)\left(\sum_{\text{j}=2}^{n}c_{j}\alpha_{\left\{ 1,j\right\} }\right)-\sum_{2\le i<j\le n}c_{i}c_{j}\alpha_{\left\{ i,j\right\} }\in\mathbb{R}\left[\Sym_{n}\right]
\]
be the difference of the two sides in the Octopus Inequality \eqref{eq:octopus}
(so $c_{i}\ge0$ for all $i=2,3,\ldots,n$). Then $\Delta|_{\D}$
is an (at most) rank-1, non-negative operator. 
\end{lem}

\begin{proof}
As in the proof of Lemma \ref{lem:The-trace-of-alpha_A in std}, we
work in the standard basis of $\D$. Then, denoting $C=\sum_{i=2}^{n}c_{i}$,
the operator $\Delta|_{\D}$ corresponds to the matrix 
\[
\left(\begin{array}{ccccc}
C^{2} & -C\cdot c_{2} & -C\cdot c_{3} & \cdots & -C\cdot c_{n}\\
-C\cdot c_{2} & c_{2}^{~2} & c_{2}c_{3} & \cdots & c_{2}c_{n}\\
-C\cdot c_{3} & c_{3}c_{2} & c_{3}^{~2} & \cdots & c_{3}c_{n}\\
\vdots & \vdots & \vdots & \ddots\\
-C\cdot c_{n} & c_{n}c_{2} & c_{n}c_{3} & \cdots & c_{n}^{~2}
\end{array}\right),
\]
which is clearly of rank at most one. The matrix is of rank-$1$ if
and only if $C\ne0$, in which case its unique nonzero (positive)
eigenvalue is given by the trace: $C^{2}+\sum_{i=2}^{n}c_{i}^{\,2}$. 
\end{proof}
We end this section with another set of representations of $\Sym_{n}$,
which are generally not irreducible, and that will be used in Section
\ref{subsec:hard cases}. For any $\mu=\left(\mu_{1},\ldots,\mu_{r}\right)\vdash n$,
we consider the subgroup $\Sym_{\mu}\defi\Sym_{\mu_{1}}\times\dotsb\times\Sym_{\mu_{r}}\subseteq\Sym_{n}$,
and let $W_{\mu}$\marginpar{$W_{\mu}$} be the induced representation
from the trivial representation of $\Sym_{\mu}$ to $\Sym_{n}$. Alternatively,
we can view $W_{\mu}$ as the representation whose elements are formal
linear combinations of cosets in $\Sym_{n}/\Sym_{\mu}$, obtained
from the action of $\Sym_{n}$ on these cosets. Note that these cosets
correspond to partitions of $\{1,\dotsc,n\}$ to disjoint sets $F_{1},\dotsc,F_{r}$
such that $|F_{i}|=\mu_{i}$ for all $i$.

The relation between the $V_{\mu}'s$ and the $W_{\mu}'s$ is given
by

\begin{equation}
W_{\mu}\cong\bigoplus_{\nu\trianglerighteq\mu}K_{\nu\mu}V_{\nu},\label{eq:Young's rule}
\end{equation}
where $\nu\trianglerighteq\mu$ is the dominance order, meaning that
$\sum_{i=1}^{j}\nu_{i}\geq\sum_{i=1}^{j}\mu_{i}$ for all $j$, and
the numbers $K_{\nu\mu}$ are positive integers known as the \emph{Kostka
numbers}. In particular, $K_{\mu\mu}=1$, namely, $V_{\mu}$ has multiplicity
one in $W_{\mu}$. The above relation is known as Young's rule \cite[Cor.~4.39]{fulton2013representation}. 

\section{The proof of Aldous' conjecture on graphs and how it may be generalized\label{sec:The-proof-of-Aldous}}

Our main results in this paper generalize both the Octopus Inequality
(Theorem \ref{thm:octopus} above) and the general framework of the
proof by Caputo-Liggett-Richthammer of Aldous' conjecture on graphs
(Theorem \ref{thm:CLR-rep-version} above). To introduce the details
of these generalizations, we first present a sketch of the proof of
the latter theorem.

\subsection{A sketch of the proof of Aldous' conjecture on graphs\label{subsec:A-sketch-of-CLR-proof}}

The proof scheme of Theorem \ref{thm:CLR} in \cite{caputo2010proof}
builds on \cite{handjani1996rate} and uses induction on the number
of vertices in the graph $G$. Assume the statement of Theorem \ref{thm:CLR}
is true for every weighted graph on $n-1$ vertices. Because the weights
on the edges of $G$ may be zero, we may assume without loss of generality
that $G$ is the full graph on $n$ vertices numbered $1,\ldots,n$
with non-negative weights $c_{i,j}\ge0$ for every $1\le i<j\le n$.
As in \eqref{eq:aldous element}, write $\Sigma=\sum_{i<j}c_{i,j}\left(i~j\right)$.
Our goal is to prove that $\lambda_{2}\left(\Sigma,\D\right)=\lambda_{2}\left(\Sigma,\R\right)$.
By Lemma \ref{lem:aldous equivalent to that std rules}, this is equivalent
to that $\lambda_{1}\left(\Sigma,V_{\mu}\right)\le\lambda_{1}\left(\Sigma,V_{\left(n-1,1\right)}\right)$
for every $\mu\vdash n$ such that $\mu\ne\left(n\right),\left(n-1,1\right)$.
Equivalently, if we let

\[
S_{G}=\sum_{i<j}c_{i,j}\alpha_{\left\{ i,j\right\} }=\sum_{i<j}c_{i,j}\left(\id-\left(i~j\right)\right)=\left(\sum_{i<j}c_{i,j}\right)-\Sigma\in\mathbb{R}\left[\Sym_{n}\right],
\]
we ought to show that 
\[
\lambda_{\min}\left(S_{G},V_{\mu}\right)\ge\lambda_{\min}\left(S_{G},V_{\left(n-1,1\right)}\right)
\]
for every $\mu\vdash n$ as above, where $\lambda_{\min}\left(S_{G},V_{\mu}\right)=\lambda_{\min}\left(S_{G}|_{V_{\mu}}\right)$
is the smallest eigenvalue of $S_{G}|_{V_{\mu}}$. Recall from Lemma
\ref{lem:alpha is non-negative} that $\alpha_{A}\ge0$ for every
$A\subseteq\left[n\right]$ so $S_{G}\ge0$. If $G$ is disconnected
(when disregarding weight-0 edges), then $\lambda_{\min}\left(S_{G},V_{\left(n-1,1\right)}\right)=0$
and we are done. So we may assume that $G$ is connected, and in particular
that for every vertex, the sum of the weights of edges incident to
it is positive. 

Define a weighted graph $H$ on the vertices $2,\ldots,n$ as follows.
For every $2\le i<j\le n$, the weight of the edge $\left\{ i,j\right\} $
is $c_{i,j}+\frac{c_{1,i}c_{1,j}}{\sum_{k=2}^{n}c_{1,k}}$. Denote
by $S_{H}$ the corresponding element of $\mathbb{R}$$\left[\Sym_{n-1}\right]$,
and by $S_{H\sqcup\left\{ 1\right\} }$ the image of $S_{H}$ through
the embedding of $\Sym_{n-1}$ in $\Sym_{n}$ as the stabilizer of
the point 1. We have 
\begin{equation}
\Delta\defi S_{G}-S_{H\sqcup\left\{ 1\right\} }=\sum_{j=2}^{n}c_{1,j}\alpha_{\left\{ 1,j\right\} }-\frac{1}{\sum_{j=2}^{n}c_{1,j}}\cdot\sum_{2\le i<j\le n}c_{1,i}c_{1,j}\alpha_{\left\{ i,j\right\} }\label{eq:G=00003DH+octopus}
\end{equation}
is (up to a constant) the ``octopus element'' from Theorem \ref{thm:octopus}.
We claim that the following equality and inequalities hold for every
$\mu\vdash n$ with $\mu\ne\left(n\right),\left(n-1,1\right)$, yielding
Theorem \ref{thm:CLR}
\begin{eqnarray}
\lambda_{\min}\left(S_{G},V_{\mu}\right) & \stackrel{\left(1\right)}{\ge} & \lambda_{\min}\left(S_{H\sqcup\left\{ 1\right\} },V_{\mu}\right)\nonumber \\
 & \stackrel{\left(2\right)}{=} & \min_{\mu'=\mu-\square}\lambda_{\min}\left(S_{H},V_{\mu'}\right)\nonumber \\
 & \stackrel{\left(3\right)}{\ge} & \lambda_{\min}\left(S_{H},V_{\left(n-2,1\right)}\right)\nonumber \\
 & \stackrel{\left(4\right)}{\ge} & \lambda_{\min}\left(S_{G},V_{\left(n-1,1\right)}\right).\label{eq:CLR proof summary}
\end{eqnarray}
Indeed,
\begin{itemize}
\item The inequality $\left(1\right)$ uses the Octopus Inequality from
Theorem \ref{thm:octopus}: $S_{G}-S_{H\sqcup\left\{ 1\right\} }\ge0$
as an operator by \eqref{eq:octopus}, which means that $\left(S_{G}-S_{H\sqcup\left\{ 1\right\} }\right)|_{V_{\mu}}\ge0$
for every $\mu\vdash n$ and $\left(1\right)$ follows.
\item The equality $\left(2\right)$ follows from the branching rule in
$\Sym_{n}$ \cite[Exer.~4.53]{fulton2013representation}: the element
$S_{H\sqcup\left\{ 1\right\} }$ is supported on the subgroup $\mathrm{Stab}_{\Sym_{n}}\left(1\right)\cong\Sym_{n-1}$
and the restriction of $V_{\mu}$ to $\Sym_{n-1}$ satisfies $V_{\mu}\cong\bigoplus_{\mu'=\mu-\square}V_{\mu'}$,
where the sum is over all Young diagrams $\mu'$ with $n-1$ squares
obtained from $\mu$ by removing a single square (which must be the
rightmost square of some row of $\mu$ and the bottom square of some
column). 
\item The inequality $\left(3\right)$ follows from the induction hypothesis
together with the observation that as $\mu\ne\left(n\right),\left(n-1,1\right)$
then after removing a single square we get that $\mu'\ne\left(n-1\right)$.
\item Finally, let us justify the inequality $\left(4\right)$. Consider
the standard representation $\D_{n-1}\cong V_{\left(n-2,1\right)}\oplus V_{\left(n-1\right)}$
of $\Sym_{n-1}$. Recall that $V_{\left(n-1\right)}=\I_{n-1}$ is
the trivial representation so $\alpha_{A}|_{V_{\left(n-1\right)}}=0$
for all $A$ and $S_{H}|_{V_{\left(n-1\right)}}=0$. As $S_{H}\ge0$,
we get that $\lambda_{\min}\left(S_{H},V_{\left(n-2,1\right)}\right)$
is the second smallest eigenvalue of $S_{H}|_{\D_{n-1}}$. Considering
$S_{H}|_{\D_{n-1}}$ and $S_{H\sqcup\left\{ 1\right\} }|_{\D_{n}}$
as matrices in the standard basis corresponding to the elements $\left\{ 2,\ldots,n\right\} $
and $\left[n\right]$, respectively, we see that $S_{H\sqcup\left\{ 1\right\} }|_{\D_{n}}$
is obtained from $S_{H}|_{\D_{n-1}}$ by adding a zero row and a zero
column, so $\lambda_{\min}\left(S_{H},V_{\left(n-2,1\right)}\right)$
is the \emph{third} smallest eigenvalue of $S_{H\sqcup\left\{ 1\right\} }|_{\D_{n}}$.
Finally, Lemma \ref{lem:Octopus is rank 1 in standard} shows that
the octopus element $\Delta$ from \eqref{eq:G=00003DH+octopus} satisfies
that $\Delta|_{\D_{n}}$ is a rank-1 (non-negative) element. By Cauchy's
interlacing formula (e.g., \cite[Lemma 3.4]{marcus2014ramanujan})
the eigenvalues $\lambda_{n}\le\dotsb\le\lambda_{1}$ of $S_{H\sqcup\left\{ 1\right\} }|_{\D_{n}}$
interlace the eigenvalues $\rho_{n}\le\dotsb\le\rho_{1}$ of $S_{G}|_{\D_{n}}=\left(S_{H\sqcup\left\{ 1\right\} }+\Delta\right)|_{\D_{n}}$
in the sense that 
\[
\lambda_{n}\le\rho_{n}\le\lambda_{n-1}\le\rho_{n-1}\le\dotsb\le\lambda_{1}\le\rho_{1},
\]
so $\lambda_{\min}\left(S_{H},V_{\left(n-2,1\right)}\right)=\lambda_{n-2}\ge\rho_{n-1}=\lambda_{\min}\left(S_{G},V_{\left(n-1,1\right)}\right)$.
\end{itemize}

\subsection{A general inequality scheme\label{subsec:A-general-inequality}}

Recall that our goal is to generalize Aldous' conjecture to hypergraphs,
in the form of Conjecture \ref{conj:Caputo's--shuffles-conjecture-Rep-version}.
This conjecture asserts that $U=\sum_{A\subseteq\left[n\right]}w_{A}J_{A}$
as in \eqref{eq:hypergraph element} satisfies $\lambda_{2}\left(U,\D\right)=\lambda_{2}\left(U,\R\right)$.
By Lemma \ref{lem:aldous equivalent to that std rules}, this is equivalent
to that $\lambda_{1}\left(U,V_{\mu}\right)\le\lambda_{1}\left(U,V_{\left(n-1,1\right)}\right)$
for all $\mu\vdash n$ with $\mu\ne\left(n\right)$,$\left(n-1,1\right)$.
As $\alpha_{A}=\left|A\right|-\frac{J_{A}}{\left(\left|A\right|-1\right)!}$,
Conjecture \ref{conj:Caputo's--shuffles-conjecture-Rep-version} is
also equivalent to that for arbitrary non-negative weights $\left\{ w_{A}\right\} _{A\subseteq\left[n\right]}$,
the element $S_{G}=\sum_{A\subseteq\left[n\right]}w_{A}\alpha_{A}$
satisfies $\lambda_{\min}\left(S_{G},V_{\mu}\right)\ge\lambda_{\min}\left(S_{G},V_{\left(n-1,1\right)}\right)$
for all $\mu\vdash n$ with $\mu\ne\left(n\right)$,$\left(n-1,1\right)$.

In an attempt to generalize to hypergraphs the argument summarized
in \eqref{eq:CLR proof summary}, one may try to use induction on
the number of vertices again. If, as above, $1$ is the ``pivotal''
vertex one tries to omit, one may hope to find an inequality similar
to the Octopus Inequality \eqref{eq:octopus}, where one side (the
larger side) is a weighted sum over all $\alpha_{A}$ with $A$ an
hyper-edge containing $1$ (and coefficients corresponding to the
weights), and the other side consists of $\alpha_{A}$'s with $A$'s
not containing 1. Such an inequality would allow us to imitate the
first step $\left(1\right)$ in \eqref{eq:CLR proof summary}. Steps
$\left(2\right)$ and $\left(3\right)$ are true by the branching
rule and induction hypothesis for general hypergraphs. Finally, in
order to recover the inequality $\left(4\right)$ in \eqref{eq:CLR proof summary},
it could be desirable that the difference between the two sides of
our inequality, when restricted to the standard representation $\D_{n}$,
be of rank 1. 

Such considerations lead us to the following proposed inequality.
Let $A_{1},\dotsc,A_{t}\subseteq\{2,3,\dotsc,n\}$, and let $c_{1},\dotsc,c_{t}\geq0$
be non-negative weights. Denote $C\defi\sum c_{i}\left|A_{i}\right|$.
When does the following inequality hold?
\begin{equation}
C\cdot\sum c_{i}\alpha_{A_{i}\cup\{1\}}\stackrel{?}{\ge}\sum(C\cdot c_{i}+c_{i}^{2})\alpha_{A_{i}}+\sum_{i<j}c_{i}c_{j}\left(\alpha_{A_{i}\cup A_{j}}+\alpha_{A_{i}\cap A_{j}}-\alpha_{A_{i}\setminus A_{j}}-\alpha_{A_{j}\setminus A_{i}}\right)\label{eq:squid}
\end{equation}
This inequality indeed generalizes the Octopus Inequality \eqref{eq:octopus}
(there $A_{1},\ldots,A_{t}$ are all singletons, and recall that $\alpha_{\left\{ i\right\} }=\alpha_{\emptyset}=0$).
It also has the property that restricted to the standard representation
$\D$, it coincides with the Octopus Inequality: this is the content
of Lemma \ref{lem:standard squid is octopus} below. Therefore, restricted
to the standard representation $\D$, the difference between the two
sides of \eqref{eq:squid} is rank-1 non-negative by Lemma \ref{lem:Octopus is rank 1 in standard},
and step $\left(4\right)$ in \eqref{eq:CLR proof summary} holds
for this generalized inequality as well.

However, there are two main problems with this attempt to generalize
the argument in \eqref{eq:CLR proof summary} to hypergraphs:
\begin{itemize}
\item First, the right hand side of \eqref{eq:squid} contains $\alpha_{A}$'s
with \emph{negative} coefficients, which means that if we simply replace
in our weighted hypergraph all the hyper-edges containing $1$ (represented
in the left hand side of \eqref{eq:squid}) with the hyper-edges as
in the right hand side of \eqref{eq:squid}, we may get a new hypergraph
with some of its weights negative, and in this regime Conjecture \ref{conj:Caputo's--shuffles-conjecture}
is not true in general. 
\item Second, and more importantly, the inequality \eqref{eq:squid} in
this general form, does always not hold. For example, it does not
hold for $n=4$, $A_{1}=\left\{ 2\right\} $, $A_{2}=\left\{ 2,3,4\right\} $
and $c_{1}=c_{2}=1$. See Section \ref{sec:Counterexamples} for more
details. 
\end{itemize}
However, the inequality \eqref{eq:squid} does hold true if certain
restrictions on the $A_{i}$'s are set. First, as mentioned above,
if all $A_{i}$'s are singletons, \eqref{eq:squid} becomes the Octopus
Inequality from Theorem \ref{thm:octopus}. More generally, \eqref{eq:squid}
holds if all the $A_{i}$'s are disjoint -- this is the content of
Theorem \ref{thm:squid - disjoint sets}; if all the $A_{i}$'s are
of size $k$ or $k-1$ and all contain some fixed subset of size $k-1$,
where $k\in\left[n\right]$ -- this is the content of Theorem \ref{thm:squid - sets of size k containing a common (k-1)-subset};
and if all sets are of size $k$ or $k-1$ and all are contained in
the some fixed subset of $\left[n\right]$ of size $k$ -- this is
the content of Theorem \ref{thm:squid - sets of cosize 1 in A0}.
\begin{rem}
\label{rem:positive coefficient on the RHS}Note that in all the cases
where we know \eqref{eq:squid} holds -- the cases covered by Theorems
\ref{thm:octopus}-\ref{thm:squid - sets of cosize 1 in A0} -- every
non-vanishing $\alpha_{A}$ on the right hand side of \eqref{eq:squid}
has a total non-negative coefficient. Indeed, in the case of disjoint
sets, $A_{i}\setminus A_{j}=A_{i}$ and the total coefficient of $\alpha_{A_{i}}$
is $c_{i}\left(c_{i}\left(\left|A_{i}\right|+1\right)+\sum_{j\ne i}c_{j}\left(\left|A_{j}\right|-1\right)\right)$
which is positive. In the remaining cases $\left|A_{i}\setminus A_{j}\right|\le1$
for all $i,j$ so $\alpha_{A_{i}\setminus A_{j}}=0$, and all remaining
summands on the right hand side of \eqref{eq:squid} have positive
coefficients. However, the inequality \eqref{eq:squid} does not necessarily
hold even if all coefficients on the right hand side are positive
--- see Example \ref{exa:squid fails even with positive coefs}.
\end{rem}

As stated above, the following lemma shows that for general sets $A_{1},\dotsc,A_{t}$
and weights $c_{1},\dotsc,c_{t}$, the inequality \eqref{eq:squid}
restricted to the standard representation becomes an Octopus Inequality.
\begin{lem}
\label{lem:standard squid is octopus}Let $A_{1},\dotsc,A_{t}\subseteq\{2,3,\dotsc,n\}$,
let $c_{1},\dotsc,c_{t}\geq0$ and denote $C=\sum c_{i}\left|A_{i}\right|$.
Consider the difference between the two sides of \eqref{eq:squid}:

\[
\Delta\defi C\cdot\sum_{i}c_{i}\alpha_{A_{i}\cup\{1\}}-\sum_{i}(C\cdot c_{i}+c_{i}^{2})\alpha_{A_{i}}-\sum_{i<j}c_{i}c_{j}\left(\alpha_{A_{i}\cup A_{j}}+\alpha_{A_{i}\cap A_{j}}-\alpha_{A_{i}\setminus A_{j}}-\alpha_{A_{j}\setminus A_{i}}\right).
\]
There exist non-negative reals $x_{2},\dotsc,x_{n}$ such that
\[
\Delta=_{\mathbf{D}}\left(\sum_{2\leq i\leq n}x_{i}\right)\left(\sum_{2\leq j\leq n}x_{j}\alpha_{\left\{ 1,j\right\} }\right)-\sum_{2\leq i<j\leq n}x_{i}x_{j}\alpha_{\left\{ i,j\right\} }.
\]
In particular, $\Delta|_{\D}$ is a non-negative rank-1 (or zero)
operator by Lemma \ref{lem:Octopus is rank 1 in standard}.
\end{lem}

\begin{proof}
For $2\leq i\leq n$ define
\[
x_{i}=\sum_{j\colon i\in A_{j}}c_{j}.
\]
Note that $\sum x_{i}=\sum_{i}\sum_{j:i\in A_{j}}c_{j}=\sum_{j}c_{j}|A_{j}|=C$.
By Lemma \ref{lem:alpha_on_standard}, we have
\begin{align*}
\sum c_{i}\left(\alpha_{A_{i}\cup\{1\}}-\alpha_{A_{i}}\right) & =_{\mathbf{D}}\frac{1}{2}\sum_{i}c_{i}\left(\sum_{j,k\in A_{i}\cup\left\{ 1\right\} }\alpha_{\left\{ j,k\right\} }-\sum_{j,k\in A_{i}}\alpha_{\left\{ j,k\right\} }\right)\\
 & =\sum_{i}c_{i}\sum_{j\in A_{i}}\alpha_{\left\{ 1,j\right\} }=\sum_{j}\alpha_{\{1,j\}}\sum_{i:j\in A_{i}}c_{i}=\sum_{j}x_{j}\alpha_{\left\{ 1,j\right\} }.
\end{align*}
It remains to prove that
\begin{equation}
\sum_{2\leq i<j\leq n}x_{i}x_{j}\alpha_{\left\{ i,j\right\} }=_{\mathbf{D}}\sum_{i}c_{i}^{2}\alpha_{A_{i}}+\sum_{i<j}c_{i}c_{j}\left(\alpha_{A_{i}\cup A_{j}}+\alpha_{A_{i}\cap A_{j}}-\alpha_{A_{i}\setminus A_{j}}-\alpha_{A_{j}\setminus A_{i}}\right).\label{eq:remains}
\end{equation}
We start by calculating
\begin{align}
\sum_{2\leq i<j\leq n}x_{i}x_{j}\alpha_{\left\{ i,j\right\} } & =\frac{1}{2}\sum_{i,j}\alpha_{\left\{ i,j\right\} }\sum_{s_{1}\colon i\in A_{s_{1}}}\sum_{s_{2}\colon j\in A_{s_{2}}}c_{s_{1}}c_{s_{2}}\nonumber \\
 & =\frac{1}{2}\sum_{s_{1},s_{2}}c_{s_{1}}c_{s_{2}}\sum_{\left(i,j\right)\in A_{s_{1}}\times A_{s_{2}}}\alpha_{\left\{ i,j\right\} }\label{eq:x_i x_j}
\end{align}
If $s_{1}=s_{2}$ we have by Lemma \ref{lem:alpha_on_standard} that
\begin{equation}
\frac{1}{2}\sum_{(i,j)\in A_{s_{1}}\times A_{s_{1}}}\alpha_{\{i,j\}}=_{\mathbf{D}}\alpha_{A_{s_{1}}}.\label{eq:s1=00003Ds2}
\end{equation}
We now deal with the case $s_{1}\ne s_{2}$. Consider the following
multiset identity, valid for any sets $A$ and $B$:
\[
A\times B+B\times A=\left(A\cup B\right)^{2}+\left(A\cap B\right)^{2}-\left(A\backslash B\right)^{2}-\left(B\backslash A\right)^{2}.
\]
By this identity and Lemma \ref{lem:alpha_on_standard}, in the case
$s_{1}\ne s_{2}$ we have that 
\begin{align}
\lefteqn{{\frac{1}{2}\left[\sum_{\left(i,j\right)\in A_{s_{1}}\times A_{s_{2}}}+\sum_{\left(i,j\right)\in A_{s_{2}}\times A_{s_{1}}}\right]\alpha_{\{i,j\}}}}\qquad\qquad\nonumber \\
 & =\frac{1}{2}\left[\sum_{(i,j)\in\left(A_{s_{1}}\cup A_{s_{2}}\right)^{2}}+\sum_{(i,j)\in\left(A_{s_{1}}\cap A_{s_{2}}\right)^{2}}-\sum_{(i,j)\in\left(A_{s_{1}}\backslash A_{s_{2}}\right)^{2}}-\sum_{(i,j)\in\left(A_{s_{1}}\backslash A_{s_{2}}\right)^{2}}\right]\alpha_{\{i,j\}}\nonumber \\
 & =_{\mathbf{D}}\alpha_{A_{s_{1}}\cup A_{s_{2}}}+\alpha_{A_{s_{1}}\cap A_{s_{2}}}-\alpha_{A_{s_{1}}\setminus A_{s_{2}}}-\alpha_{A_{s_{2}}\setminus A_{s_{1}}}\label{eq:s1 ne s2}
\end{align}
Putting \eqref{eq:s1=00003Ds2} and \eqref{eq:s1 ne s2} into \eqref{eq:x_i x_j},
we get \eqref{eq:remains} as needed.
\end{proof}

\subsection{Reducing Theorem \ref{thm:triangles and the like - Rep version}
to Theorem \ref{thm:squid - sets of size k containing a common (k-1)-subset}}

We can now prove our main result about hypergraphs --- Theorem \ref{thm:triangles and the like - Rep version}
--- assuming Theorem \ref{thm:squid - sets of size k containing a common (k-1)-subset},
which generalizes the Octopus Inequality. Theorem \ref{thm:squid - sets of size k containing a common (k-1)-subset}
is proved in Section \ref{sec:sets with large intersection}.
\begin{proof}[Proof of Theorem \ref{thm:triangles and the like - Rep version} assuming
Theorem \ref{thm:squid - sets of size k containing a common (k-1)-subset}]
 Let $B\subseteq\left[n\right]$. We need to show that 

\[
U=\sum_{\substack{B\subseteq A\subseteq\left[n\right]\colon\\
\left|A\setminus B\right|\le2
}
}w_{A}J_{A}\in\mathbb{R}\left[\Sym_{n}\right]
\]
satisfies $\lambda_{2}\left(U,\D\right)=\lambda_{2}\left(U,\R\right)$
for arbitrary non-negative weights $\left\{ w_{A}\right\} _{A}$.
By Lemma \ref{lem:aldous equivalent to that std rules} and the definition
of $\alpha_{A}$, this is equivalent to that for arbitrary non-negative
weights $\left\{ w_{A}\right\} _{A}$, the element 

\[
S_{G}=\sum_{\substack{B\subseteq A\subseteq\left[n\right]\colon\\
\left|A\setminus B\right|\le2
}
}w_{A}\alpha_{A}\in\mathbb{R}\left[\Sym_{n}\right]
\]
satisfies that for every $\mu\vdash n$ with $\mu\ne\left(n\right),\left(n-1,1\right)$
\[
\lambda_{\min}\left(S_{G},V_{\mu}\right)\ge\lambda_{\min}\left(S_{G},V_{\left(n-1,1\right)}\right).
\]
We prove this by induction on $n-\left|B\right|$. If $B=\left[n\right]$,
then $S_{G}=w_{B}\alpha_{B}$, by Lemma \ref{lem:alpha_A in S_A}
$\lambda_{\min}\left(S_{G},V_{\mu}\right)=\omega_{B}n$ for every
$\mu\vdash n$ with $\mu\ne\left(n\right)$ and we are done.

Now assume that $\left|B\right|<n$ and assume the theorem is true
whenever $n-\left|B\right|$ is smaller. As in Section \ref{subsec:A-sketch-of-CLR-proof},
if the hypergraph $G$ is disconnected (disregarding weight-0 hyper-edges),
then $S_{G}\ge0$ and $\lambda_{\min}\left(S_{G},V_{\left(n-1,1\right)}\right)=0$
and we are done. So we may assume that $G$ is connected. In particular,
every vertex is contained in a hyper-edge with a positive weight.
Let $v\in\left[n\right]\setminus B$ be a vertex outside $B$. Without
loss of generality, we may assume that $v=1$. 

Denote $t=n-1-\left|B\right|$, $A_{0}=B$ and $c_{0}=w_{B\cup\left\{ 1\right\} }$.
Let $\left[n\right]\setminus B=\left\{ 1,m_{1},\ldots,m_{t}\right\} $
and for every $i=1,\ldots,t$ let $A_{i}=B\cup\left\{ m_{i}\right\} $
and $c_{i}=w_{A_{i}\cup\left\{ 1\right\} }$. Also denote $C=\sum_{i=0}^{t}c_{i}\left|A_{i}\right|$,
and note that $C>0$ by our assumption that every vertex is contained
in a positive-weight hyper-edge. 

Now define a hypergraph $H$ on the vertices $\left\{ 2,\ldots,n\right\} $
with hyper-edges consisting of all subsets $A$ containing $B$ and
satisfying $\left|A\setminus B\right|\le2$, with the following non-negative
weights $\overline{w}_{A}$:
\begin{equation}
\overline{w}_{A}=\begin{cases}
w_{A}+c_{0}+\frac{1}{C}\left(c_{0}^{~2}+\sum_{0\le i<j\le t}c_{i}c_{j}\right) & \mathrm{if}~A=B,\\
w_{A}+c_{i}+\frac{1}{C}\left(c_{i}^{~2}+c_{0}c_{i}\right) & \mathrm{if}~A=B\cup\left\{ m_{i}\right\} ,1\le i\le t\\
w_{A}+\frac{1}{C}c_{i}c_{j} & \mathrm{if}~A=B\cup\left\{ m_{i},m_{j}\right\} ,1\le i<j\le t.
\end{cases}\label{eq:new_weights}
\end{equation}
Note that the hypergraph $H$ satisfies the assumptions of Theorem
\ref{thm:triangles and the like - Rep version} with vertex set $\left\{ 2,\ldots,n\right\} $
and the same special subset $B$. We denote 
\[
S_{H\cup\left\{ 1\right\} }=\sum_{\substack{B\subseteq A\subseteq\left\{ 2,\ldots,n\right\} \colon\\
\left|A\setminus B\right|\le2
}
}\overline{w}_{A}\alpha_{A}\in\mathbb{R}\left[\Sym_{n}\right],
\]
and $S_{H}$ the corresponding element of $\mathbb{R}\left[\Sym_{n-1}\right]$
(here $\Sym_{n-1}$ is embedded in $\Sym_{n}$ as the stabilizer of
the point $1$). By the induction hypothesis, $\lambda_{\min}\left(S_{H},V_{\mu'}\right)\ge\lambda_{\min}\left(S_{H},V_{\left(n-2,1\right)}\right)$
for every $\mu'\vdash n-1$ with $\mu'\ne\left(n-1\right),\left(n-2,1\right)$.
The proof can now be completed as in \eqref{eq:CLR proof summary}.
Indeed, for every $\mu\vdash n$ we have
\begin{align}
\lambda_{\min}\left(S_{G},V_{\mu}\right) & \stackrel{\left(1\right)}{\ge}\lambda_{\min}\left(S_{H\sqcup\left\{ 1\right\} },V_{\mu}\right)\stackrel{\left(2\right)}{=}\min_{\mu'=\mu-\square}\lambda_{\min}\left(S_{H},V_{\mu'}\right)\nonumber \\
 & \stackrel{\left(3\right)}{\ge}\lambda_{\min}\left(S_{H},V_{\left(n-2,1\right)}\right)\stackrel{\left(4\right)}{\ge}\lambda_{\min}\left(S_{G},V_{\left(n-1,1\right)}\right).\label{eq:proof-of-triangles}
\end{align}
The inequality $\left(1\right)$ follows from that $S_{G}-S_{H\cup\left\{ 1\right\} }$
is, up to a multiplicative constant, the operator which is non-negative
by Theorem \ref{thm:squid - sets of size k containing a common (k-1)-subset}.
The equality $\left(2\right)$ follows from the branching rule in
$\Sym_{n}$ \cite[Exer.~4.53]{fulton2013representation}. The inequality
$\left(3\right)$ holds by the induction hypothesis. Finally, $S_{G}-S_{H\cup\left\{ 1\right\} }$
is a special case of the difference between the two sides of \eqref{eq:squid},
so $\left(S_{G}-S_{H\cup\left\{ 1\right\} }\right)|_{\D}$ coincides
with an octopus element by Lemma \ref{lem:standard squid is octopus},
and thus the inequality $\left(4\right)$ is true similarly to the
inequality $\left(4\right)$ in \eqref{eq:CLR proof summary}.
\end{proof}

\section{Proof of Theorem \ref{thm:squid - disjoint sets}\label{sec:Proof-of- disjoint sets}}

The first generalization of the Octopus Inequality that we prove is
Theorem \ref{thm:squid - disjoint sets} dealing with disjoint sets.
We will need the following lemma. Recall that $\ker\alpha_{A}$ denotes
the kernel of the operator $\alpha_{A}$ in $\mathbb{R}\left[\Sym_{n}\right]$.
\begin{lem}
\label{lem:alpha A+1 in ker alpha_A}For every $A\subseteq\left\{ 2,\ldots,n\right\} $,
the following equality of operators holds inside $\ker\alpha_{A}$:

\begin{equation}
\alpha_{A\cup\left\{ 1\right\} }=\sum_{a\in A}\alpha_{\left\{ 1,a\right\} }.\label{eq:alpha_A+1 inside ker alpha_A}
\end{equation}
\end{lem}

\begin{proof}
The case $A=\emptyset$ is trivial, so assume otherwise. Consider
the operator 
\[
T\defi\left(\alpha_{A\cup\left\{ 1\right\} }-\sum_{\left\{ x,y\right\} \in\binom{A\cup\left\{ 1\right\} }{2}}\alpha_{\left\{ x,y\right\} }\right)\left(\alpha_{A}-\left|A\right|\right)\in\mathbb{R}\left[\Sym_{n}\right].
\]
Note that $T$ is supported on $A\cup\left\{ 1\right\} $, and we
claim it vanishes in $S_{A\cup\left\{ 1\right\} }$. Indeed, set $k=\left|A\right|+1$.
In the standard representation $\D_{k}=V_{\left(k\right)}\oplus V_{\left(k-1,1\right)}$
Lemma \ref{lem:alpha_on_standard} yields that $\alpha_{A\cup\left\{ 1\right\} }-\sum_{\left\{ x,y\right\} \in\binom{A\cup\left\{ 1\right\} }{2}}\alpha_{\left\{ x,y\right\} }=0$.
For any $\mu\vdash k$ with $\mu\ne\left(k\right),\left(k-1,1\right)$
we have $\alpha_{A}-\left|A\right|=0$ by Corollary \ref{cor:alpha is constant when many squares outside first row}.
By Lemma \ref{lem:positive_on_Sk} we conclude that $T\equiv0$ on
$\mathbb{R}\left[\Sym_{n}\right]$.

Restricted to $\ker\alpha_{A}\subseteq\mathbb{R}\left[\Sym_{n}\right]$,
the operator $\alpha_{A}-\left|A\right|$ acts by multiplication by
a non-zero scalar, so $\alpha_{A\cup\left\{ 1\right\} }=\sum_{\left\{ x,y\right\} \in\binom{A\cup\left\{ 1\right\} }{2}}\alpha_{\left\{ x,y\right\} }$.
Finally, by Lemma \ref{lem:if A subset of B}, $\ker\alpha_{A}\subseteq\ker\alpha_{\left\{ x,y\right\} }$
whenever $x,y\in A$, namely, inside $\ker\alpha_{A}$ we have $\alpha_{\left\{ x,y\right\} }=0$.
We conclude that inside $\ker\alpha_{A}$
\[
\alpha_{A\cup\left\{ 1\right\} }=\sum_{\left\{ x,y\right\} \in\binom{A\cup\left\{ 1\right\} }{2}}\alpha_{\left\{ x,y\right\} }=\sum_{a\in A}\alpha_{\left\{ 1,a\right\} }.\qedhere
\]
\end{proof}
Recall that in the statement of Theorem \ref{thm:squid - disjoint sets},
$A_{1},\ldots,A_{t}\subseteq\left\{ 2,\ldots,n\right\} $ are pairwise
disjoint subsets, $c_{1},\ldots,c_{t}$ are non-negative numbers,
and the theorem states the following inequality of operators in $\mathbb{R}\left[\Sym_{n}\right]$:
\begin{equation}
\left(\sum_{i=1}^{t}c_{i}\left|A_{i}\right|\right)\left(\sum_{j=1}^{t}c_{j}\left(\alpha_{A_{j}\cup\left\{ 1\right\} }-\alpha_{A_{j}}\right)\right)\ge\sum_{i=1}^{t}c_{i}^{~2}\alpha_{A_{i}}+\sum_{1\le i<j\le t}c_{i}c_{j}\left(\alpha_{A_{i}\cup A_{j}}-\alpha_{A_{i}}-\alpha_{A_{j}}\right).\label{eq:again squid for disjoint sets}
\end{equation}

\begin{proof}[Proof of Theorem \ref{thm:squid - disjoint sets}]
 Denote $C=\sum_{i=1}^{t}c_{i}\left|A_{i}\right|$. Rearranging \eqref{eq:again squid for disjoint sets},
we need to prove that 
\begin{equation}
\sum_{i=1}^{t}c_{i}\left(C\cdot\alpha_{A_{i}\cup\left\{ 1\right\} }-\left(C+c_{i}-\sum_{j\ne i}c_{j}\right)\alpha_{A_{i}}\right)\ge\sum_{1\le i<j\le t}c_{i}c_{j}\alpha_{A_{i}\cup A_{j}}.\label{eq:rearranging}
\end{equation}
Denote by $L$ the left hand side of \eqref{eq:rearranging} and by
$R$ the right hand side. A key observation is that for every $i=1,\ldots,t$,
the operator $\alpha_{A_{i}}$ commutes with $\alpha_{B}$ for every
$\alpha_{B}$ appearing in \eqref{eq:rearranging}: indeed, every
such $B$ either \emph{contains} $A_{i}$ and then $\alpha_{A_{i}}\alpha_{B}=\alpha_{B}\alpha_{A_{i}}$
by Lemma \ref{lem:if A subset of B}\eqref{enu:product alpha_a alpha_b},
or \emph{is disjoint} from $A_{i}$ and then $\alpha_{A_{i}}\alpha_{B}=\alpha_{B}\alpha_{A_{i}}$
because $\sigma\tau=\tau\sigma$ for every $\sigma,\tau\in\Sym_{n}$
with $\mathrm{supp}\left(\sigma\right)\subseteq A_{i}$, $\mathrm{supp}\left(\tau\right)\subseteq B$.
In particular, $\alpha_{A_{1}},\ldots,\alpha_{A_{t}}$ can be diagonalized
simultaneously, and we may split the space $\mathbb{R}\left[\Sym_{n}\right]$
to $2^{t}$ subspaces
\[
\mathbb{R}\left[\Sym_{n}\right]=\bigoplus_{z\subseteq\left[t\right]}V_{z},
\]
where for every $z\subseteq\left[t\right]$, 
\[
V_{z}=\left(\bigcap_{i\in z}\ker\alpha_{A_{i}}\right)\cap\left(\bigcap_{i\notin z}\ker\alpha_{A_{i}}^{~\perp}\right).
\]
Note that $\ker\alpha_{A_{i}}^{~~\perp}$ is nothing more than the
eigenspace of $\alpha_{A_{i}}$ with eigenvalue $\left|A_{i}\right|$,
by Lemma \ref{lem:if A subset of B}\eqref{enu:eigenvalues of alpha_a}.
Moreover, as $L$ and $R$ commute with $A_{i}$ for every $i$, they
preserve the eigenspaces of $A_{i}$, and thus also preserve $V_{z}$
for each $z\subseteq\left[t\right]$. It is therefore enough to prove
the inequality \eqref{eq:rearranging} on $V_{z}$ for every $z$
separately.

Fix, therefore, some $z\subseteq\left[t\right]$. By Lemma \ref{lem:if A subset of B}\eqref{enu:kernels},
for every $i\notin z$, $\ker\alpha_{A_{i}\cup\left\{ 1\right\} }\subseteq\ker\alpha_{A_{i}}$
so $\ker\alpha_{A_{i}}^{~\perp}\subseteq\ker\alpha_{A_{i}\cup\left\{ 1\right\} }^{~\perp}$
(recall that $\alpha_{A_{i}}$ and $\alpha_{A_{i}\cup\left\{ 1\right\} }$
commute), and by Lemma \ref{lem:if A subset of B}\eqref{enu:eigenvalues of alpha_a},
inside $\ker\alpha_{A_{i}}^{~\perp}$ we have $\alpha_{A_{i}}\equiv\left|A_{i}\right|$,
$\alpha_{A_{i}\cup\left\{ 1\right\} }\equiv\left|A_{i}\right|+1$
and similarly $\alpha_{A_{i}\cup A_{j}}\equiv\left|A_{i}\right|+\left|A_{j}\right|$.
Restricted to $V_{z}$, the $i^{th}$ summand in the left hand side
of \eqref{eq:rearranging} is thus (still for $i\not\in z$)
\begin{align*}
 & c_{i}\left[C\cdot\alpha_{A_{i}\cup\left\{ 1\right\} }-\left(C+c_{i}-\sum_{j\ne i}c_{j}\right)\alpha_{A_{i}}\right]\\
 & =c_{i}\left(C\cdot\left(\left|A_{i}\right|+1\right)-\left(C+c_{i}-\sum_{j\ne i}c_{j}\right)\left|A_{i}\right|\right)\\
 & =c_{i}\left(C-c_{i}\left|A_{i}\right|+\sum_{j\ne i}c_{j}\left|A_{i}\right|\right)=c_{i}\left(\sum_{j}c_{j}\left|A_{j}\right|-c_{i}\left|A_{i}\right|+\sum_{j\ne i}c_{j}\left|A_{i}\right|\right)\\
 & =c_{i}\sum_{j\ne i}c_{j}\left(\left|A_{j}\right|+\left|A_{i}\right|\right)=\sum_{j\ne i}c_{i}c_{j}\alpha_{A_{i}\cup A_{j}}.
\end{align*}
So restricting to $V_{z}$, if we subtract this last equality from
\eqref{eq:rearranging} for every $i\notin z$, we are left to prove
that
\[
\sum_{i\in z}c_{i}\left(C\cdot\alpha_{A_{i}\cup\left\{ 1\right\} }-\left(C+c_{i}-\sum_{j\ne i}c_{j}\right)\alpha_{A_{i}}\right)\ge\sum_{\left\{ i,j\right\} \in\binom{z}{2}}c_{i}c_{j}\alpha_{A_{i}\cup A_{j}}-\sum_{\left\{ i,j\right\} \in\binom{\left[t\right]\setminus z}{2}}c_{i}c_{j}\alpha_{A_{i}\cup A_{j}}.
\]
Using the non-negativity of $c_{i}c_{j}\alpha_{A_{i}\cup A_{j}}$
and that $\alpha_{A_{i}}\equiv0$ in $\ker\alpha_{A_{i}}$, it is
enough to prove that 
\begin{equation}
C\cdot\sum_{i\in z}c_{i}\alpha_{A_{i}\cup\left\{ 1\right\} }\ge\sum_{\left\{ i,j\right\} \in\binom{z}{2}}c_{i}c_{j}\alpha_{A_{i}\cup A_{j}}.\label{eq:enough to prove}
\end{equation}
If $z=\emptyset$ we are done. Otherwise, we conclude by the following
sequence of equalities and inequalities in $V_{z}$. 
\begin{eqnarray*}
C\cdot\sum_{i\in z}c_{i}\alpha_{A_{i}\cup\left\{ 1\right\} } & \stackrel{\left(1\right)}{=} & C\cdot\sum_{i\in z}c_{i}\sum_{a\in A_{i}}\alpha_{\left\{ 1,a\right\} }\\
 & \stackrel{\left(2\right)}{\ge} & \frac{C}{\sum_{i\in z}c_{i}\left|A_{i}\right|}\cdot\left[\sum_{i\in z}c_{i}^{2}\sum_{\left\{ x,y\right\} \in\binom{A_{i}}{2}}\alpha_{\left\{ x,y\right\} }+\sum_{\left\{ i,j\right\} \in\binom{z}{2}}c_{i}c_{j}\sum_{\substack{x\in A_{i}\\
y\in A_{j}
}
}\alpha_{\left\{ x,y\right\} }\right]\\
 & \stackrel{\left(3\right)}{\ge} & \sum_{\left\{ i,j\right\} \in\binom{z}{2}}c_{i}c_{j}\sum_{\substack{x\in A_{i}\\
y\in A_{j}
}
}\alpha_{\left\{ x,y\right\} }\\
 & \stackrel{\left(4\right)}{=} & \sum_{\left\{ i,j\right\} \in\binom{z}{2}}c_{i}c_{j}\sum_{\left\{ x,y\right\} \in\binom{A_{i}\cup A_{j}}{2}}\alpha_{\left\{ x,y\right\} }\\
 & \stackrel{\left(5\right)}{\ge} & \sum_{\left\{ i,j\right\} \in\binom{z}{2}}c_{i}c_{j}\alpha_{A_{i}\cup A_{j}}.
\end{eqnarray*}
Equality $\left(1\right)$ holds inside $V_{z}$ by Lemma \ref{lem:alpha A+1 in ker alpha_A}.
Inequality $\left(2\right)$ is an instance of the Octopus Inequality
from Theorem \ref{thm:octopus}. Inequality $\left(3\right)$ follows
from that $0<\sum_{i\in z}c_{i}\left|A_{i}\right|\le\sum_{i=1}^{t}c_{i}\left|A_{i}\right|=C$
and that if $x,y\in A_{i}$ then $\alpha_{\left\{ x,y\right\} }\equiv0$
inside $\ker\alpha_{A_{i}}$. The latter argument also gives equality
$\left(4\right)$, and inequality $\left(5\right)$ follows from inequality
\eqref{eq:alpha_A le sum of alpha of pairs} in Lemma \ref{lem:alpha_on_standard}.
\end{proof}

\section{Proof of Theorem \ref{thm:squid - sets of size k containing a common (k-1)-subset}
\label{sec:sets with large intersection}}

In this section we prove Theorem \ref{thm:squid - sets of size k containing a common (k-1)-subset}
--- our generalization of the Octopus Inequality for sets with large
intersections, which is used in the proof of our main result on hypergraphs
(Theorem \ref{thm:triangles and the like - Rep version}). Recall
that the setting is that $A_{0},A_{1},\ldots,A_{t}\subseteq\left\{ 2,\ldots,n\right\} $
are distinct subsets with $A_{i}\supseteq A_{0}$ and $\left|A_{i}\setminus A_{0}\right|=1$
for all $i=1,\ldots,t$, and that $c_{0},\ldots,c_{t}$ are non-negative
real numbers. We denote the size of $A_{1},\ldots,A_{t}$ by $k$,
so that $\left|A_{0}\right|=k-1$ and $k\in\mathbb{Z}_{\ge1}$. The
case $k=1$ is the setting of the original octopus inequality, so
we may assume that $k\ge2$.

We rearrange (\ref{eq:squid - sets with large intersection}) and
get that we need to prove 
\begin{equation}
U\ge0,\label{eq:squid_for_k}
\end{equation}
where
\begin{equation}
\begin{aligned}U\defi & \sum_{i=0}^{t}c_{i}^{2}\left(\left|A_{i}\right|\alpha_{A_{i}\cup\{1\}}-\left(\left|A_{i}\right|+1\right)\alpha_{A_{i}}\right)+\\
 & +\sum_{0\leq i<j\le t}c_{i}c_{j}\left(\left|A_{i}\right|\left(\alpha_{A_{j}\cup\{1\}}-\alpha_{A_{j}}\right)+\left|A_{j}\right|\left(\alpha_{A_{i}\cup\{1\}}-\alpha_{A_{i}}\right)-\alpha_{A_{i}\cup A_{j}}-\alpha_{A_{0}}\right).
\end{aligned}
\label{eq:def of U}
\end{equation}
For the proof of \eqref{eq:squid_for_k} we will use the squaring
strategy of \cite{caputo2010proof}. This strategy is based on the
following idea: by Lemma \ref{lem:standard squid is octopus}, the
restriction of $U$ to the standard representation is equal to the
difference between the sides of an octopus inequality. By Lemma \ref{lem:Octopus is rank 1 in standard},
we conclude that $U|_{\D}$ is of rank at most 1. Let $d=\text{tr}(U|_{\mathbf{D}})$.
Then $d$ is the only (possibly) non-zero eigenvalue of $U|_{\mathbf{D}}$.
In addition, by the Octopus Inequality \eqref{eq:octopus}, $U|_{\D}\geq0$,
so $d\geq0$.

We will prove not only that $U\ge0$, but also that $U\le d$. In
other words, we will prove the operator inequality $\left(U-d/2\right)^{2}\leq\left(d/2\right)^{2}$,
or, equivalently 
\begin{equation}
\left(d/2\right)^{2}-\left(U-d/2\right)^{2}\geq0,\label{eq:main inequality}
\end{equation}
from which \eqref{eq:squid_for_k} immediately follows. Let us start
by evaluating $d$. By Lemma \ref{lem:The-trace-of-alpha_A in std},
$\tr\left(\alpha_{A}|_{\mathbf{D}}\right)=|A|\left(|A|-1\right)$,
and recall that $\left|A_{0}\right|=k-1$ and $\left|A_{i}\right|=k$
for $i=1,\ldots,t$. So,
\begin{align*}
d & =\tr(U|_{\mathbf{D}})=\sum_{i=0}^{t}c_{i}^{~2}\left(\left|A_{i}\right|^{2}\left(\left|A_{i}\right|+1\right)-\left(\left|A_{i}\right|+1\right)\left|A_{i}\right|\left(\left|A_{i}\right|-1\right)\right)\\
 & +\sum_{0\le i<j\le t}c_{i}c_{j}\left(4\left|A_{i}\right|\left|A_{j}\right|-\left|A_{i}\cup A_{j}\right|\left(\left|A_{i}\cup A_{j}\right|-1\right)-\left(k-1\right)\left(k-2\right)\right)\\
 & =c_{0}^{~2}k\left(k-1\right)+\sum_{i=1}^{t}c_{i}^{~2}k\left(k+1\right)+\sum_{j=1}^{t}c_{0}c_{j}\left(2k^{2}-2\right)+\sum_{1\le i<j\le t}c_{i}c_{j}\left(2k^{2}+2k-2\right).
\end{align*}
We will prove \eqref{eq:main inequality} by expanding both squares
and collecting terms by monomials in the variables $\{c_{i}\}.$ Let
us write $U=\sum_{0\leq i\le j\le t}c_{i}c_{j}U_{ij}$, where
\begin{equation}
U_{ij}=\begin{cases}
\left(k-1\right)\alpha_{A_{0}\cup\left\{ 1\right\} }-k\alpha_{A_{0}} & 0=i=j\\
k\alpha_{A_{i}\cup\{1\}}-\left(k+1\right)\alpha_{A_{i}} & 1\le i=j\\
k\alpha_{A_{0}\cup\{1\}}+\left(k-1\right)\alpha_{A_{j}\cup\{1\}}-\left(k+1\right)\alpha_{A_{0}}-k\alpha_{A_{j}} & 0=i<j\\
k\left(\alpha_{A_{i}\cup\{1\}}+\alpha_{A_{j}\cup\{1\}}-\alpha_{A_{i}}-\alpha_{A_{j}}\right)-\alpha_{A_{i}\cup A_{j}}-\alpha_{A_{0}} & 1\le i<j
\end{cases}\label{eq:Uij}
\end{equation}
We also have $\frac{d}{2}=\sum_{0\leq i\le j\le t}c_{i}c_{j}w_{ij}$
with
\begin{equation}
w_{ij}=\begin{cases}
\frac{\left(k-1\right)k}{2} & 0=i=j\\
\frac{k\left(k+1\right)}{2} & 1\le i=j\\
k^{2}-1 & 0=i<j\\
k^{2}+k-1 & 1\le i<j
\end{cases}\label{eq:w_ij}
\end{equation}
Now \eqref{eq:main inequality} becomes
\begin{equation}
\left[\sum_{0\leq i\le j\le t}c_{i}c_{j}w_{ij}\right]^{2}-\left[\sum_{0\leq i\le j\le t}c_{i}c_{j}\left(U_{ij}-w_{ij}\right)\right]^{2}\ge0.\label{eq:elaborate_inequality}
\end{equation}
Let us write the left hand side as
\[
\sum_{0\leq i\le j\le\ell\le m\le t}c_{i}c_{j}c_{\ell}c_{m}\Gamma_{ij\ell m}.
\]
We will show that $\Gamma_{ij\ell m}\geq0$ for all $i,j,\ell,m$,
and \eqref{eq:main inequality} will follow.

For convenience, we will denote $\Gamma_{ij\ell m}=\Gamma_{i'j'\ell'm'}$
whenever $(i,j,\ell,m)$ (any quadruple of non-negative integers)
is a rearrangement of $(i',j',\ell',m')$ and $i'\leq j'\leq\ell'\leq m'$,
and similarly for $U_{ij}$ and for $w_{ij}$.

When expanding the squares in \eqref{eq:elaborate_inequality}, we
have to be mindful of the order of multiplication. To that end, we
define
\begin{defn}
\label{def:david-star}For operators $A$ and $B$ denote by $A\david B$
the expression $\frac{AB+BA}{2}$. This operation is called the \textbf{Jordan
product}.
\end{defn}

We get

\begin{equation}
\Gamma_{ij\ell m}=\sum_{\substack{e\leq f,e'\leq f'\colon\\
\{e,f\}\cup\{e',f'\}=\left\{ i,j,\ell,m\right\} 
}
}\Lambda_{efe'f'}\label{eq:gamma_by_lambda}
\end{equation}
where the union $\{e,f\}\cup\{e',f'\}=\left\{ i,j,\ell,m\right\} $
is a union of multisets, and

\begin{equation}
\Lambda_{efe'f'}=w_{ef}w_{e'f'}-\left(U_{ef}-w_{ef}\right)\david\left(U_{e'f'}-w_{e'f'}\right).\label{eq:def_of_lambda}
\end{equation}
For example, if $i<j<\ell<m$ then $\Gamma_{ij\ell m}=\Lambda_{ij\ell m}+\Lambda_{\ell mij}+\Lambda_{i\ell jm}+\Lambda_{jmi\ell}+\Lambda_{imj\ell}+\Lambda_{j\ell im}$,
and if $i<j$ then $\Gamma_{iijj}=\Lambda_{iijj}+\Lambda_{jjii}+\Lambda_{ijij}$.

Fixing $i,j,\ell,m$, denote \marginpar{$E$} 
\[
E=A_{i}\cup A_{j}\cup A_{\ell}\cup A_{m}\cup\left\{ 1\right\} 
\]
and note that $k\leq|E|\leq k+4$ and that the group ring element
$\Gamma_{ij\ell m}$ is supported on $\Sym_{E}$. By Lemma \ref{lem:positive_on_Sk},
it is enough to show that $\Gamma_{ij\ell m}\geq0$ in the group ring
of $\Sym_{E}$. Therefore, we need to show that $\Gamma_{ij\ell m}\geq_{V_{\mu}}0$
for any $\mu\vdash|E|$. In the following subsections we will prove
this claim under various assumptions on $\mu_{1}$, the length of
the first row of $\mu$: in Section \ref{subsec:The-easy-cases} we
deal with the most straightforward cases where $\mu_{1}\le k-1$ or
$\left|E\right|-1\le\mu_{1}$; in Section \ref{subsec:lambda_1=00003Dk}
we prove the inequality when $\mu_{1}=k$; and in Section \ref{subsec:hard cases}
we deal with the remaining cases, namely, $k+1\le\mu_{1}\le\left|E\right|-2$.

The following lemma will be helpful in all cases.
\begin{lem}
\label{lem:Ker_alpha_is_enough}For every $i,j,\ell,m$, the subspaces
$\ker\alpha_{A_{0}},\ker\alpha_{A_{0}}^{~\perp}\subseteq\mathbb{R}\left[\Sym_{E}\right]$
are invariant under the operator $\Gamma_{\ijlm}$, and $\Gamma_{\ijlm}|_{\ker\alpha_{A_{0}}^{~\perp}}=0$.
Thus, in order to show that $\Gamma_{ij\ell m}$ is non-negative,
it is enough to prove that $\Gamma_{\ijlm}|_{\ker\alpha_{A_{0}}}$
is non-negative.
\end{lem}

\begin{proof}
Recall from Lemma \ref{lem:if A subset of B}\eqref{enu:eigenvalues of alpha_a}
that the eigenvalues of $\alpha_{A}$ are $0$ and $\left|A\right|$,
so $\mathbb{R}\left[\Sym_{E}\right]=\ker\alpha_{A_{0}}\oplus\ker\alpha_{A_{0}}^{~\perp}$,
where $\ker\alpha_{A_{0}}^{~\perp}=\ker\left(\alpha_{A_{0}}-\left(k-1\right)\right)$.
Moreover, as $A_{0}\subseteq B$ for every $B\subseteq E$ such that
$\alpha_{B}$ appears in $\Gamma_{\ijlm}$, then $\Gamma_{\ijlm}$
commutes with $\alpha_{A_{0}}$ by Lemma \ref{lem:if A subset of B}\eqref{enu:product alpha_a alpha_b},
so $\Gamma_{\ijlm}$ preserves the eigenspaces of $\alpha_{A_{0}}$.
So it is enough to prove that $\Gamma_{\ijlm}|_{\ker\alpha_{A_{0}}}\ge0$
and that $\Gamma_{\ijlm}|_{\ker\alpha_{A_{0}}^{~\perp}}\ge0$.

Now, for every $B\subseteq E$ such that $\alpha_{B}$ appears in
the expression defining $\Gamma_{\ijlm}$, as $\text{\ensuremath{\ker\alpha_{B}\subseteq\ker\alpha_{A_{0}}} }$
and $\alpha_{A_{0}}$ and $\alpha_{B}$ commute by Lemma \ref{lem:if A subset of B},
we also have $\ker\alpha_{A_{0}}^{~\perp}\subseteq\ker\alpha_{B}^{~\perp}=\ker\left(\alpha_{B}-\left|B\right|\right)$.
So inside $\ker\alpha_{A_{0}}^{~\perp}$ we have $\alpha_{B}\equiv\left|B\right|$
for every $B\supseteq A_{0}$. Substituting $\alpha_{B}$ with $\left|B\right|$
in \eqref{eq:Uij}, we see that for all $e,f\in\left\{ i,j,\ell,m\right\} $
we have $U_{ef}|_{\ker\alpha_{A_{0}}^{~\perp}}=0$, so $U|_{\ker\alpha_{A_{0}}^{~\perp}}=0$
and thus $\Gamma_{\ijlm}|_{\ker\alpha_{A_{0}}^{~\perp}}=0$.
\end{proof}

\subsection{The easy cases: $\mu_{1}\le k-1$ or $\mu_{1}\ge\left|E\right|-1$\label{subsec:The-easy-cases}}
\begin{lem}
\label{Lem:1st row<k-1}An equality $\Gamma_{ij\ell m}=_{V_{\mu}}0$
holds for all $i,j,\ell,m$ and every $\mu\vdash\left|E\right|$ whenever
$\mu_{1}\le k-2$.
\end{lem}

\begin{proof}
If $\left|A_{0}\right|=k-1\le1$ there are no such partitions $\mu$.
So assume that $\left|A_{0}\right|\ge2$. By Corollary \ref{cor:alpha is constant when many squares outside first row},
for every such $\mu$ we have $\alpha_{A_{0}}=_{V_{\mu}}|A_{0}|$,
so $V_{\mu}\subseteq\ker\alpha_{A_{0}}^{~\perp}$ and we are done
by Lemma \ref{lem:Ker_alpha_is_enough}.
\end{proof}
\begin{lem}
\label{Lem:1st row=00003Dk-1}The inequality $\Gamma_{ij\ell m}\geq_{V_{\mu}}0$
holds for all $i,j,\ell,m$ and every $\mu\vdash\left|E\right|$ with
$\mu_{1}=k-1$.
\end{lem}

\begin{proof}
By Lemma \ref{lem:Ker_alpha_is_enough}, it is enough to prove that
$\Gamma_{ij\ell m}$ is non-negative on $\ker\alpha_{A_{0}}\subseteq V_{\mu}$.
In that subspace $\alpha_{A_{0}}=0$, and by Corollary \ref{cor:alpha is constant when many squares outside first row},
$\alpha_{B}=\left|B\right|$ for every $A_{0}\subsetneqq B\subseteq E$.
Hence all the operators $\alpha_{A_{i}}$ are scalars, and consequently,
so are the operators $U_{ef}$ for all $e,f\in\left\{ i,j,\ell,m\right\} $.
By \eqref{eq:Uij} we have

\[
U_{ef}=\begin{cases}
k\left(k-1\right) & 0=e=f\\
0 & 1\le e=f\\
k^{2}-1 & 0=e<f\\
k-1 & 1\le e<f
\end{cases}.
\]
Comparing these constants with $w_{ef}$ from \eqref{eq:w_ij}, we
see that $0\leq U_{ef}\leq2w_{ef}$ for all $e,f\in\left\{ i,j,\ell,m\right\} $,
so $|U_{ef}-w_{ef}|\leq w_{ef}$. We conclude that

\[
\left|(U_{ef}-w_{ef})\david(U_{e'f'}-w_{e'f'})\right|=\left|(U_{ef}-w_{ef})\cdot(U_{e'f'}-w_{e'f'})\right|\leq w_{ef}w_{e'f'},
\]
so 
\[
\Lambda_{efe'f'}=w_{ef}w_{e'f'}-\left(U_{ef}-w_{ef}\right)\david\left(U_{e'f'}-w_{e'f'}\right)\geq0.
\]
By \eqref{eq:gamma_by_lambda} we have $\Gamma_{\ijlm}\geq_{V_{\mu}}0$.
\end{proof}
\begin{lem}
\label{lem:mu_1 =00005Cge |E|-1}The equality $\Gamma_{ij\ell m}=_{V_{\mu}}0$
holds for all $i,j,\ell,m$ and every $\mu\vdash\left|E\right|$ with
$\mu_{1}\ge\left|E\right|-1$. 
\end{lem}

\begin{proof}
As mentioned above, $U$ is of rank $1$ on $\mathbf{D=\D}_{n}$.
Hence, the only non-zero eigenvalue of $U$ in $\mathbf{D}$ is $\tr\left(U|_{\mathbf{D}}\right)=d$.
Therefore, $(U-d/2)^{2}$ has only $(d/2)^{2}$ as eigenvalue, so
$\left(d/2\right)^{2}-\left(U-d/2\right)^{2}=_{\mathbf{D}}0$, and
$\sum_{i\le j\le\ell\le m}c_{i}c_{j}c_{\ell}c_{m}\Gamma_{ij\ell m}=_{\mathbf{D}}0$
for any choice of non-negative values for $\left\{ c_{i}\right\} $.
A polynomial in the $\left\{ c_{i}\right\} $ that vanishes on every
non-negative values of the $\left\{ c_{i}\right\} $ must be the zero
polynomial, so $\Gamma_{\ijlm}=_{\D}0$. Finally, recall that $\mathbf{D}\cong V_{\left(n-1,1\right)}\oplus V_{\left(n\right)},$so
$\Gamma_{\ijlm}=_{V_{\left(n-1,1\right)}}0$ and $\Gamma_{\ijlm}=_{V_{\left(n\right)}}0$.
By the branching rule, we get $\Gamma_{\ijlm}=_{V_{\mu}}0$.
\end{proof}
At this point we have established that $\Gamma_{\ijlm}\ge0$ when
$\left|E\right|\le k+1$, namely, for $\Gamma_{0000}$, $\Gamma_{000i}$,
$\Gamma_{00ii}$, $\Gamma_{0iii}$ and $\Gamma_{iiii}$ for every
$1\le i\le t$.

\subsection{The case $\mu_{1}=k$\label{subsec:lambda_1=00003Dk}}

In this subsection we deal with the case that $\mu\vdash\left|E\right|$
with $\mu_{1}=k$. We assume that $k<\left|E\right|-1$, for otherwise
the case $\mu_{1}=k$ is contained in Lemma \ref{lem:mu_1 =00005Cge |E|-1}.
Lemma \ref{lem:mu1=00003Dk i,j,l,m positive} below proves that $\Gamma_{\ijlm}\ge_{V_{\mu}}0$
when $i,j,\ell,m$ are all positive and not necessarily distinct,
Lemma \ref{lem:mu1=00003Dk 00ij} addresses the case of $\Gamma_{00ij}$
with $i$ and $j$ distinct, Lemma \ref{lem:mu1=00003Dk 0iij } the
case of $\Gamma_{0iij}$ with $i$ and $j$ distinct, and Lemma \ref{lem:mu1=00003Dk 0ijl}
the case of $\Gamma_{0ij\ell}$ with $i,j,\ell$ distinct. As mentioned
above, the remaining cases -- $\Gamma_{0000}$, $\Gamma_{000i}$,
$\Gamma_{00ii}$, $\Gamma_{0iii}$ and $\Gamma_{iiii}$ -- have already
been established in Section \ref{subsec:The-easy-cases}.
\begin{lem}
\label{lem:mu1=00003Dk i,j,l,m positive}Assume that $i,j,\ell,m>0$
and $\mu\vdash\left|E\right|$ with $\mu_{1}=k$. Then $\Gamma_{ij\ell m}\ge_{V_{\mu}}0$.
\end{lem}

\begin{proof}
By Lemma \ref{lem:Ker_alpha_is_enough}, it is enough to prove that
$\Gamma_{ij\ell m}$ is non-negative on $\ker\alpha_{A_{0}}\subseteq V_{\mu}$.
By Corollary \ref{cor:alpha is constant when many squares outside first row},
we have $\alpha_{B}=_{V_{\mu}}k+1$ for every $B\subset E$ with $\left|B\right|=k+1$.
As we now show, this fact guarantees that 
\begin{equation}
\left\Vert U_{ef}-w_{ef}\right\Vert \le w_{ef}\label{eq:uij-wij has norm less than wij}
\end{equation}
for all $e,f>0$. Recall that by Lemma \ref{lem:if A subset of B},
$\left\Vert \alpha_{A}-\frac{|A|}{2}\right\Vert =\text{\ensuremath{\frac{|A|}{2}}}$
for any $A\subseteq E$.

Let us consider first the case $1\le e=f$. We have by \eqref{eq:Uij}
and \eqref{eq:w_ij} that 
\begin{eqnarray*}
\left\Vert U_{ee}-w_{ee}\right\Vert  & = & \left\Vert k\alpha_{A_{e}\cup\{1\}}-(k+1)\alpha_{A_{e}}-\frac{k\left(k+1\right)}{2}\right\Vert \\
 & = & \left\Vert k\left(k+1\right)-\left(k+1\right)\alpha_{A_{e}}-\frac{k\left(k+1\right)}{2}\right\Vert \\
 & = & \left\Vert \left(k+1\right)\left(\frac{k}{2}-\alpha_{A_{e}}\right)\right\Vert =\left(k+1\right)\frac{k}{2}=w_{ee}.
\end{eqnarray*}
In the case $1\le e<f$:
\begin{eqnarray*}
\left\Vert U_{ef}-w_{ef}\right\Vert  & = & \left\Vert k(\alpha_{A_{e}\cup\{1\}}+\alpha_{A_{f}\cup\{1\}}-\alpha_{A_{e}}-\alpha_{A_{f}})-\alpha_{A_{e}\cup A_{f}}-\alpha_{A_{0}}-\left(k^{2}+k-1\right)\right\Vert \\
 & = & \left\Vert k\left(2k+2-\alpha_{A_{e}}-\alpha_{A_{f}}\right)-\left(k+1\right)-\left(k^{2}+k-1\right)\right\Vert \\
 & = & \left\Vert k\left(k-\alpha_{A_{e}}-\alpha_{A_{f}}\right)\right\Vert =k\left\Vert \left(\frac{k}{2}-\alpha_{A_{e}}\right)+\left(\frac{k}{2}-\alpha_{A_{f}}\right)\right\Vert \\
 & \leq & k\left\Vert \frac{k}{2}-\alpha_{A_{e}}\right\Vert +k\left\Vert \frac{k}{2}-\alpha_{A_{f}}\right\Vert =k^{2}\le k^{2}+k-1=w_{ef}.
\end{eqnarray*}
As in the proof of Lemma \ref{Lem:1st row=00003Dk-1}, we conclude
that
\[
\left\Vert \left(U_{ef}-w_{ef}\right)\david\left(U_{e'f'}-w_{e'f'}\right)\right\Vert \leq w_{ef}w_{e'f'},
\]
so $\Lambda_{efe'f'}\geq_{V_{\mu}}0$ whenever $\left\{ e,f,e',f'\right\} =\left\{ i,j,\ell,m\right\} $,
and by \eqref{eq:gamma_by_lambda}, $\Gamma_{ijlm}\geq_{V_{\mu}}0$.
\end{proof}
\begin{lem}
\label{lem:mu1=00003Dk 00ij}Assume that $i,j>0$ are distinct and
that $\mu\vdash\left|E\right|$ with $\mu_{1}=k$. Then $\Gamma_{00ij}\ge_{V_{\mu}}0$.
\end{lem}

\begin{proof}
By Corollary \ref{cor:alpha is constant when many squares outside first row},
$\alpha_{A_{i}\cup\left\{ 1\right\} }=_{V_{\mu}}\alpha_{A_{j}\cup\left\{ 1\right\} }=_{V_{\mu}}\alpha_{A_{i}\cup A_{j}}=_{V_{\mu}}k+1$.
As before, it is enough to prove the non-negativity on $V_{\mu}\cap\ker\alpha_{A_{0}}$,
so we assume that $\alpha_{A_{0}}=0$. We therefore get from \eqref{eq:Uij}
\begin{eqnarray*}
U_{00} & = & \left(k-1\right)\alpha_{A_{0}\cup\left\{ 1\right\} },\\
U_{0i} & = & k\left(\alpha_{A_{0}\cup\left\{ 1\right\} }-\alpha_{A_{i}}\right)+k^{2}-1,\\
U_{ij} & = & 2k^{2}+k-1-k\left(\alpha_{A_{i}}+\alpha_{A_{j}}\right).
\end{eqnarray*}
With \eqref{eq:w_ij} we have
\begin{eqnarray*}
U_{00}-w_{00} & = & \left(k-1\right)\left(\alpha_{A_{0}\cup\left\{ 1\right\} }-\frac{k}{2}\right),\\
U_{0i}-w_{0i} & = & k\left(\alpha_{A_{0}\cup\left\{ 1\right\} }-\alpha_{A_{i}}\right),\\
U_{ij}-w_{ij} & = & k^{2}-k\left(\alpha_{A_{i}}+\alpha_{A_{j}}\right).
\end{eqnarray*}
By Lemma \ref{lem:if A subset of B}\eqref{enu:norm-of-alpha_a} we
get the operator norm bounds
\begin{eqnarray*}
\left\Vert U_{00}-w_{00}\right\Vert  & = & \frac{\left(k-1\right)k}{2},\\
\left\Vert U_{0i}-w_{0i}\right\Vert  & = & k\left\Vert \left(\alpha_{A_{0}\cup\left\{ 1\right\} }-\frac{k}{2}\right)-\left(\alpha_{A_{i}}-\frac{k}{2}\right)\right\Vert \le k^{2},\\
\left\Vert U_{ij}-w_{ij}\right\Vert  & = & k\left\Vert \left(\frac{k}{2}-\alpha_{A_{i}}\right)+\left(\frac{k}{2}-\alpha_{A_{j}}\right)\right\Vert \le k^{2}.
\end{eqnarray*}
By \eqref{eq:gamma_by_lambda} we have that in $\ker\alpha_{A_{0}}\cap V_{\mu}$
\begin{eqnarray*}
\Gamma_{00ij} & = & 2(\Lambda_{00ij}+\Lambda_{0i0j})\\
 & = & 2w_{00}w_{ij}-2\left(U_{00}-w_{00}\right)\david\left(U_{ij}-w_{ij}\right)+\\
 &  & 2w_{0i}w_{0j}-2\left(U_{0i}-w_{0i}\right)\david\left(U_{0j}-w_{0j}\right),
\end{eqnarray*}
so it is enough to show that 
\[
\left\Vert \left(U_{00}-w_{00}\right)\david\left(U_{ij}-w_{ij}\right)+\left(U_{0i}-w_{0i}\right)\david\left(U_{0j}-w_{0j}\right)\right\Vert \le w_{00}w_{ij}+w_{0i}w_{0j}.
\]
The left hand side is bounded from above by $\frac{\left(k-1\right)k}{2}\cdot k^{2}+k^{4}$.
The right hand side is equal to $\frac{\left(k-1\right)k}{2}\cdot\left(k^{2}+k-1\right)+\left(k^{2}-1\right)^{2}$.
So our claim holds if
\[
\frac{\left(k-1\right)k}{2}\cdot k^{2}+k^{4}\le\frac{\left(k-1\right)k}{2}\cdot\left(k^{2}+k-1\right)+\left(k^{2}-1\right)^{2},
\]
which is equivalent to
\[
k^{3}-6k^{2}+k+2\ge0.
\]
This inequality holds for $k\ge6$. For\footnote{Recall that that we assume throughout the proof of Theorem \ref{thm:squid - sets of size k containing a common (k-1)-subset}
that $k\ge2$. In fact, when $k=1$, a direct computation gives that
$U_{0i}=w_{0i}=0$ for arbitrary $i$, so $\Gamma_{0ij\ell}=0$ for
arbitrary $i,j,\ell$.} $k=2,3,4,5$ we verified by computer that $\Gamma_{00ij}\ge_{V_{\mu}}0$
-- see Appendix \ref{sec:All-remaining-cases small k}. 
\end{proof}
\begin{lem}
\label{lem:mu1=00003Dk 0iij }Assume that $i,j>0$ are distinct and
that $\mu\vdash\left|E\right|$ with $\mu_{1}=k$. Then $\Gamma_{0iij}\ge_{V_{\mu}}0$.
\end{lem}

\begin{proof}
This proof is similar to the proof of Lemma \ref{lem:mu1=00003Dk 00ij}.
We have $\left|E\right|=k+2$. As before, it is enough to prove the
non-negativity of $\Gamma_{0iij}$ on $V_{\mu}\cap\ker\alpha_{A_{0}}$.
In this subspace, $\alpha_{A_{i}\cup\left\{ 1\right\} }=_{V_{\mu}}\alpha_{A_{j}\cup\left\{ 1\right\} }=_{V_{\mu}}\alpha_{A_{i}\cup A_{j}}=_{V_{\mu}}k+1$
and $\alpha_{A_{0}}=_{V_{\mu}}0$.

We have seen in the proof of Lemma \ref{lem:mu1=00003Dk 00ij}, that
in this situation we have $\left\Vert U_{0i}-w_{0i}\right\Vert \leq k^{2}$
and $\left\Vert U_{ij}-w_{ij}\right\Vert \leq k^{2}$. Let us estimate
$\left\Vert U_{ii}-w_{ii}\right\Vert $. By \eqref{eq:Uij},

\[
U_{ii}=k\alpha_{A_{i}\cup\left\{ 1\right\} }-\left(k+1\right)\alpha_{A_{i}}=k\left(k+1\right)-\left(k+1\right)\alpha_{A_{i}},
\]
and with \eqref{eq:w_ij} we have
\[
U_{ii}-w_{ii}=\frac{k\left(k+1\right)}{2}-\left(k+1\right)\alpha_{A_{i}}.
\]
As $\left\Vert \alpha_{A}-\frac{\left|A\right|}{2}\right\Vert =\frac{\left|A\right|}{2}$
by Lemma \ref{lem:if A subset of B}\eqref{enu:norm-of-alpha_a} we
get
\[
\left\Vert U_{ii}-w_{ii}\right\Vert =(k+1)\left\Vert \frac{k}{2}-\alpha_{A_{i}}\right\Vert \le\frac{k\left(k+1\right)}{2}.
\]
By \eqref{eq:gamma_by_lambda} (assuming without loss of generality
that $i<j$ -- this only affects the notation),
\begin{eqnarray*}
\Gamma_{0iij} & = & 2(\Lambda_{0iij}+\Lambda_{0jii})\\
 & = & 2w_{0i}w_{ij}-2\left(U_{0i}-w_{0i}\right)\david\left(U_{ij}-w_{ij}\right)+\\
 &  & 2w_{0j}w_{ii}-2\left(U_{0j}-w_{0j}\right)\david\left(U_{ii}-w_{ii}\right),
\end{eqnarray*}
so it is enough to show that 
\begin{equation}
\left\Vert \left(U_{0i}-w_{0i}\right)\david\left(U_{ij}-w_{ij}\right)+\left(U_{0j}-w_{0j}\right)\david\left(U_{ii}-w_{ii}\right)\right\Vert \le w_{0i}w_{ij}+w_{0j}w_{ii}.\label{eq:bound on norm}
\end{equation}
The left hand side is bounded from above by $k^{4}+k^{2}\cdot\frac{k(k+1)}{2}$.
The right hand side is equal to $\left(k^{2}-1\right)\cdot\left(k^{2}+k-1\right)+\left(k^{2}-1\right)\cdot\frac{k(k+1)}{2}$.
So our claim holds if
\[
k^{4}+\frac{k^{3}\left(k+1\right)}{2}\le\left(k^{2}-1\right)\cdot\left(k^{2}+k-1\right)+\left(k^{2}-1\right)\cdot\frac{k\left(k+1\right)}{2},
\]
which is equivalent to
\[
2k^{3}-5k^{2}-3k+2\ge0.
\]
This inequality holds for $k\ge3$. For $k=2$ we verified by computer
that $\Gamma_{0iij}\ge0$ -- see Appendix \ref{sec:All-remaining-cases small k}.
\end{proof}
\begin{lem}
\label{lem:mu1=00003Dk 0ijl}Assume that $i,j,\ell>0$ are distinct,
and that $\mu\vdash\left|E\right|$ with $\mu_{1}=k$. Then $\Gamma_{0ij\ell}\ge_{V_{\mu}}0$.
\end{lem}

\begin{proof}
This proof is similar to the proof of the last two lemmas. As before,
we consider only the restriction of $\Gamma_{0ij\ell}$ to $\ker\alpha_{A_{0}}\cap V_{\mu}$,
where $\alpha_{A_{i}\cup\left\{ 1\right\} }=_{V_{\mu}}\alpha_{A_{j}\cup\left\{ 1\right\} }=_{V_{\mu}}\alpha_{A_{i}\cup A_{j}}=_{V_{\mu}}k+1$
and $\alpha_{A_{0}}=_{V_{\mu}}0$.

By \eqref{eq:gamma_by_lambda},
\begin{eqnarray*}
\Gamma_{0ij\ell} & = & 2(\Lambda_{0ij\ell}+\Lambda_{0ji\ell}+\Lambda_{0\ell ij})\\
 & = & 2w_{0i}w_{j\ell}-2\left(U_{0i}-w_{0i}\right)\david\left(U_{j\ell}-w_{j\ell}\right)+\\
 &  & 2w_{0j}w_{i\ell}-2\left(U_{0j}-w_{0j}\right)\david\left(U_{i\ell}-w_{i\ell}\right)+\\
 &  & 2w_{0\ell}w_{ij}-2\left(U_{0\ell}-w_{0\ell}\right)\david\left(U_{ij}-w_{ij}\right),
\end{eqnarray*}
so it is enough to show that
\[
\begin{array}{l}
\left\Vert \left(U_{0i}-w_{0i}\right)\david\left(U_{j\ell}-w_{j\ell}\right)+\left(U_{0j}-w_{0j}\right)\david\left(U_{i\ell}-w_{i\ell}\right)+\left(U_{0\ell}-w_{0\ell}\right)\david\left(U_{ij}-w_{ij}\right)\right\Vert \\
\le w_{0i}w_{j\ell}+w_{0j}w_{i\ell}+w_{0\ell}w_{ij}.
\end{array}
\]

We have seen in the proof of Lemma \ref{lem:mu1=00003Dk 00ij}, that
in this situation we have $\left\Vert U_{0e}-w_{0e}\right\Vert \leq k^{2}$
and $\left\Vert U_{ef}-w_{ef}\right\Vert \leq k^{2}$ for all distinct
$e,f\in\{i,j,\ell\}$. Hence, the left hand side is bounded from above
by $3k^{4}$. The right hand side is equal to $3\cdot\left(k^{2}-1\right)\cdot\left(k^{2}+k-1\right)$.
So our claim holds if
\[
k^{4}\le\left(k^{2}-1\right)\cdot\left(k^{2}+k-1\right),
\]
which is equivalent to
\[
k^{3}-2k^{2}-k+1\ge0.
\]
This inequality holds for $k\ge3$. For $k=2$ we verified by computer
that in these representations $\Gamma_{0ij\ell}\ge0$ -- see Appendix
\ref{sec:All-remaining-cases small k}. 
\end{proof}
At this point we have established that $\Gamma_{\ijlm}\ge0$ also
when $\left|E\right|=k+2$, namely, for $\Gamma_{00ij}$, $\Gamma_{0iij}$,
$\Gamma_{iiij}$, and $\Gamma_{iijj}$ for every $i\ne j$ with $i,j\in\left\{ 1,\ldots,t\right\} $.

\subsection{The case $k+1\le\mu_{1}\le\left|E\right|-2$\label{subsec:hard cases}}

Recall that $\left|E\right|\le k+4$, so if $k+1\le\mu_{1}\le\left|E\right|-2$
then necessarily $\left|E\right|=k+3$ or $\left|E\right|=k+4$. Namely,
the cases we ought to consider are the following three:
\begin{enumerate}
\item Prove that $\Gamma_{ij\ell m}\ge_{V_{\mu}}0$ when $i,j,\ell,m>0$
are distinct and $\mu\vdash k+4$ with $\mu_{1}=k+1$ or $\mu_{1}=k+2$.
\item Prove that $\Gamma_{iij\ell}\ge_{V_{\mu}}0$ when $i,j,\ell>0$ are
distinct and $\mu\vdash k+3$ with $\mu_{1}=k+1$.
\item Prove that $\Gamma_{0ij\ell}\ge_{V_{\mu}}0$ when $i,j,\ell>0$ are
distinct and $\mu\vdash k+3$ with $\mu_{1}=k+1$.
\end{enumerate}
Without loss of generality, we assume that $A_{0}=\left\{ 2,\ldots,k\right\} $
and $A_{i}=A_{0}\cup\left\{ k+1\right\} $, $A_{j}=A_{0}\cup\left\{ k+2\right\} $,
$A_{\ell}=A_{0}\cup\left\{ k+3\right\} $ and $A_{m}=A_{0}\cup\left\{ k+4\right\} $.

By Young's rule \eqref{eq:Young's rule}, all the representations
in the first claim\footnote{For $k\ge2,$ the irreducible representations in the first claim are
$V_{\mu}$ for $\mu=\left(k+1,3\right)$, $\mu=\left(k+1,2,1\right)$,
$\mu=\left(k+1,1,1,1\right)$, $\mu=\left(k+2,2\right)$ and $\mu=\left(k+2,1,1\right)$.} are subrepresentations of the reducible representation $W_{\left(k+1,1,1,1\right)}$
-- the action of $\Sym_{k+4}$ on ordered triplets of distinct numbers
from $\left\{ 1,\ldots,k+4\right\} $. Likewise, we need to prove
that $\Gamma_{iij\ell}\ge0$ and $\Gamma_{0ij\ell}\ge0$ on irreducible
representations $V_{\mu}$ of $\Sym_{k+3}$ with $\mu_{1}=k+1$, all
of which are subrepresentations of the reducible representation $W_{\left(k+1,1,1\right)}$
-- the action of $\Sym_{k+3}$ on ordered pairs of distinct numbers
from $\left\{ 1,\ldots,k+3\right\} $.

Our strategy for all three cases is similar, and we illustrate it
by the first one -- $\Gamma_{ij\ell m}$. By Lemma \ref{lem:Ker_alpha_is_enough},
$V_{k}:=W_{\left(k+1,1,1,1\right)}\cap\ker\alpha_{A_{0}}$ is invariant
under $\Gamma_{ij\ell m}$, and it is enough to prove that $\Gamma_{ij\ell m}\ge0$
in this $V_{k}$. Since $\frac{1}{k-1}\alpha_{A_{0}}=\mathrm{id}-\frac{1}{\left(k-1\right)!}\sum_{\sigma\in\mathrm{Sym}\left(A_{0}\right)}\sigma$,
$V_{k}$ consists of linear combinations of triplets which are invariant
under $\Sym_{A_{0}}$. Therefore, for all $k\ge4$, the following
list of vectors forms a linear basis for $V_{k}$:
\begin{itemize}
\item $v_{x,y,z}\defi\left(x,y,z\right)$ for distinct $x,y,z\in\left\{ 1,k+1,k+2,k+3,k+4\right\} $.
\item ${\displaystyle v_{x,y,*}\defi\sum_{s\in\left\{ 2,\ldots,k\right\} }\left(x,y,s\right)}$
and similarly $v_{x,*,y}$ and $v_{*,x,y}$ for distinct $x,y\in\{1,k+1,k+2,k+3,\linebreak[0]k+4\}$.
\item ${\displaystyle v_{x,*,*}\defi\sum_{\substack{s,t\in\left\{ 2,\ldots,k\right\} \colon\\
s\ne t
}
}\left(x,s,t\right)}$ and similarly $v_{*,x,*}$ and $v_{*,*,x}$ for $x\in\{1,k+1,k+2,k+3,k+4\}$.
\item ${\displaystyle v_{*,*,*}\defi\sum_{s,t,r\in\left\{ 2,\ldots,k\right\} \colon~s,t,r~\mathrm{distinct}}\left(s,t,r\right)}$.
\end{itemize}
In particular, the dimension of $V_{k}$ is the same, $136$, for
every $k\ge4$. (For $k\le3$ some of these elements vanish. For example,
for $k=3$, the element $V_{*,*,*}=0$. We check these cases separately
and directly by computer --- see Appendix \ref{sec:All-remaining-cases small k}.)

Our next goal is to construct an explicit $136\times136$ matrix giving
the action of $\Gamma_{ij\ell m}$ on $V_{k}$ with respect to this
given basis. We will show that the entries of this matrix are rational
functions of the parameter $k$, namely, elements in $\mathbb{Q}\left(k\right)$.
As $\Gamma_{ij\ell m}\,|\,_{V_{k}}$ is a (quadratic) polynomial in
the operators $\alpha_{B}\,|\,_{V_{k}}$ with $A_{0}\subseteq B\subseteq\{1,\ldots,k+4\}$,
it is enough to describe the $136\times136$ matrix corresponding
to each such $\alpha_{B}$.

Fix a subset $B$ as above. Let $v=v_{x,y,z}$ be a basis element,
where $x,y,z\in\{1,k+1,k+2,k+3,k+4,*\}$. Pick a triplet $\left(a,b,c\right)\in\{1,\dotsc,k+4\}^{3}$
in the support of the sum defining $v$ and define $O_{B}(v)$ as
the $\Sym_{B}$-orbit of $(a,b,c)$. Since $A_{0}\subseteq B,$ we
have that $O_{B}(v)$ does not depend on the choice of $(a,b,c)$
but only on $B$ and $v$. Let $T_{B}\left(v\right)$ be the set of
all basis elements $v'$ supported on $O_{B}\left(v\right)$. So $O_{B}\left(v\right)\subseteq\left[k+4\right]^{3}$
and $T_{B}\left(v\right)$ is a subset of the 136 basis elements. 
\begin{lem}
\label{lem:action of alpha_B}Fix $B$ and $v=v_{x,y,z}$ a basis
element as above. Let $n_{*}\left(v\right)$ denote the number of
stars among $x,y,z$ and $n_{B}\left(v\right)$ denote the number
of $x,y,z$ which are either a star or an element of $B\setminus A_{0}$.
The action of $\alpha_{B}$ on $v$ is given by 
\begin{equation}
\alpha_{B}\cdot v=\left|B\right|v-\frac{(\left|B\right|-n_{B}(v))!\cdot\left(k-1\right)\cdots\left(k-n_{*}\left(v\right)\right)}{(|B|-1)!}\sum_{v'\in T_{B}\left(v\right)}v'.\label{eq:alpha_B * v}
\end{equation}
\end{lem}

\begin{proof}
Recall that $J_{B}=\sum_{\sigma\in\mathrm{Sym}(B)}\sigma$, so $\alpha_{B}=\left|B\right|-\frac{J_{B}}{\left(\left|B\right|-1\right)!}$.
It is enough to prove that 
\begin{equation}
\frac{J_{B}}{\left(\left|B\right|-1\right)!}\cdot v=\frac{(\left|B\right|-n_{B}(v))!\cdot\left(k-1\right)\cdots\left(k-n_{*}\left(v\right)\right)}{(\left|B\right|-1)!}\sum_{v'\in T_{B}\left(v\right)}v'.\label{eq:J_B.v}
\end{equation}
We have $|O_{B}(v)|=\left|B\right|!/(|B|-n_{B}(v))!$, since within
the triplet defining $v$, only the $n_{B}(v)$ elements belonging
to $B$ can be moved by $\Sym_{B}$, and they all move to distinct
places within $B$. Hence, for every $t\in O_{B}\left(v\right)$,
\[
\left(\sum_{\sigma\in\Sym_{B}}\sigma\right)t=\frac{\left|\Sym_{B}\right|}{\left|O_{B}\left(v\right)\right|}\left(\sum_{s\in O_{B}\left(v\right)}s\right)=\frac{\left|B\right|!}{\left|B\right|!/(|B|-n_{B}(v))!}\left(\sum_{s\in O_{B}\left(v\right)}s\right).
\]
Summing the above equality over all $\left(k-1\right)\cdots\left(k-n_{*}\left(v\right)\right)$
triplets in the support of the sum defining $v$, and as $\sum_{s\in O_{B}\left(v\right)}s=\sum_{v'\in T_{B}\left(v\right)}v'$,
we obtain \eqref{eq:J_B.v}.
\end{proof}
Denote by $M\left(\Gamma_{ij\ell m}\right)$ the matrix of the operator
$\Gamma_{ij\ell m}$ in its action on $V_{k}$, with respect to the
above basis.
\begin{lem}
The entries of $M\left(\Gamma_{ij\ell m}\right)$ are rational functions
of $k$, and $k^{2}M(\Gamma_{ijlm})\in M_{136}(\mathbb{Z}[k])$.
\end{lem}

\begin{proof}
By the formulas \eqref{eq:gamma_by_lambda}, \eqref{eq:Uij} and \eqref{eq:w_ij},
$\Gamma_{ij\ell m}$ is a quadratic polynomial in the operators $\alpha_{B}$
for some sets \textbf{$A_{0}=\{2,\dotsc,k\}\subseteq B\subseteq\{1,\dotsc,k+4\}$
}such that $k-1\leq|B|\leq k+1$, with coefficients in $\mathbb{Z}[k]$.
For each such $B$, consider the coefficient 
\[
\frac{\left(|B|-n_{B}(v)\right)!\left(k-1\right)\cdots\left(k-n_{*}\left(v\right)\right)}{\left(|B|-1\right)!}=\frac{|B|\left(k-1\right)\cdots\left(k-n_{*}\left(v\right)\right)}{\left|B\right|\left(|B|-1\right)\cdots\left(|B|-n_{B}\left(v\right)+1\right)}
\]
in \eqref{eq:alpha_B * v}. We have $|B|=k-1+|B\backslash A_{0}|$
and $n_{B}(v)\leq n_{*}(v)+|B\backslash A_{0}|,$ hence 
\[
|B|-n_{B}(v)+1\geq k-1+|B\backslash A_{0}|-\left(n_{*}(v)+|B\backslash A_{0}|\right)+1=k-n_{*}(v).
\]
As $|B|\leq k+1$, all the terms in the denominator, with the possible
exception of $k$ (if $\left|B\right|-1=k$), are canceled by terms
in the numerator. We conclude that the matrix of $\alpha_{B}$ is
in $\frac{1}{k}M_{136}(\mathbb{Z}[k])$, hence $M(\Gamma_{ij\ell m})\in\frac{1}{k^{2}}M_{136}(\mathbb{Z}[k])$.
\end{proof}
Using a computer (see Appendix \ref{sec:code hard cases for almost all k}),
we calculated the matrix $A=k^{2}M(\Gamma_{ijlm})$, and found its
characteristic polynomial $p(k,t)\in\mathbb{Z}[k][t]$. We then calculated
$q(k,t):=p(k+3,-t)$. As it turns out, all the coefficients of $q$
are non-negative. Therefore, for each $k\geq1$, $q(k,t)$ has no
positive real roots, hence for each $k\geq4$, the characteristic
polynomial of $M(\Gamma_{ij\ell m})$ has no negative roots, so the
operator $\Gamma_{ij\ell m}$ has no negative eigenvalues. We conclude
that $\Gamma_{ij\ell m}\geq_{V_{k}}0$ for $k\ge4$. The cases $k=2,3$
were verified separately and directly -- see Appendix \ref{sec:All-remaining-cases small k}.
\begin{rem}
Our computer computations are all carried over the integers, so there
are no issues of numerical mistakes.
\end{rem}

The other two cases -- of $\Gamma_{iij\ell}$ and of $\Gamma_{0ij\ell}$
-- are proven similarly, with the computerized proof for $k\ge3$
appearing in Appendix \ref{subsec:code for Gamma iijl and 0ijl},
and the computerized proof for $k=2$ appearing in Appendix \ref{sec:All-remaining-cases small k}.

This concludes the proof of Theorem \ref{thm:squid - sets of size k containing a common (k-1)-subset}.

\section{Proof of Theorem \ref{thm:squid - sets of cosize 1 in A0} \label{sec:n-1 sets}}

Recall that in the statement of this theorem there is a set $A_{0}\subseteq\left\{ 2,\ldots,n\right\} $,
there are distinct sets $A_{1},\ldots,A_{t}$ satisfying $A_{i}\subseteq A_{0}$
and $\left|A_{0}\setminus A_{i}\right|=1$ for $i=1,\ldots,t$, and
there are non-negative numbers $c_{0},c_{1},\ldots,c_{t}\ge0$. We
ought to prove the following inequality of operators in $\mathbb{R}\left[\Sym_{n}\right]$:
\begin{equation}
\left(\sum_{i=0}^{t}c_{i}\left|A_{i}\right|\right)\left(\sum_{j=0}^{t}c_{j}\left(\alpha_{A_{j}\cup\left\{ 1\right\} }-\alpha_{A_{j}}\right)\right)\ge\sum_{i=0}^{t}c_{i}^{~2}\alpha_{A_{i}}+\sum_{0\le i<j\le t}c_{i}c_{j}\left(\alpha_{A_{0}}+\alpha_{A_{i}\cap A_{j}}\right).\label{eq:squid - sets of cosize 1 - again}
\end{equation}

\begin{proof}[Proof of Theorem \ref{thm:squid - sets of cosize 1 in A0}]
 Note that all the elements in \eqref{eq:squid - sets of cosize 1 - again}
are supported on $A_{0}\cup\left\{ 1\right\} $, so by Lemma \ref{lem:positive_on_Sk},
it is enough to prove that the inequality \eqref{eq:squid - sets of cosize 1 - again}
is satisfied in $\mathbb{R}\left[\Sym_{A_{0}\cup\left\{ 1\right\} }\right]$.
Without loss of generality, assume that $A_{0}=\left\{ 2,\ldots,n\right\} $.
Denote 
\[
C=\sum_{i=0}^{t}c_{i}\left|A_{i}\right|=\left(n-1\right)c_{0}+\left(n-2\right)\sum_{i=1}^{t}c_{i},
\]
and define $S\defi\sum_{i=1}^{t}c_{i}$, $T\defi\sum_{1\le i<j\le t}c_{i}c_{j}$
and $Q\defi\sum_{i=1}^{t}c_{i}^{2}$, so $C=\left(n-1\right)c_{0}+\left(n-2\right)S$
and $S^{2}=2T+Q$. Now rewrite \eqref{eq:squid - sets of cosize 1 - again}
as 
\begin{multline}
C\cdot\left(\sum_{i=0}^{t}c_{i}\alpha_{A_{i}\cup\left\{ 1\right\} }\right)-\left(C\cdot c_{0}+c_{0}^{2}+\sum_{0\le i<j\le t}c_{i}c_{j}\right)\alpha_{A_{0}}\\
\ge\left(\sum_{i=1}^{t}c_{i}\left(C+c_{i}+c_{0}\right)\alpha_{A_{i}}\right)+\sum_{1\le i<j\le t}c_{i}c_{j}\alpha_{A_{i}\cap A_{j}},\label{eq:constants on the left}
\end{multline}
so that all the $\alpha_{B}$'s with $\left|B\right|\ge n-1$ are
on the left hand side, and those with $\left|B\right|\le n-2$ are
on the right hand side. By Lemma \ref{lem:standard squid is octopus},
the inequality we ought to prove is satisfied in $\D$, namely in
$V_{\left(n\right)}$ and $V_{\left(n-1,1\right)}$. It remains to
prove it for $V_{\mu}$ where $\mu\vdash n$ and $\mu_{1}\le n-2$.
Under these assumptions, by Corollary \ref{cor:alpha is constant when many squares outside first row},
if $\left|B\right|\ge n-1$ then $\alpha_{B}=_{V_{\mu}}\left|B\right|$,
so the left hand side of \eqref{eq:constants on the left} is equal
to 
\begin{align*}
\mathrm{LHS} & =C\left(c_{0}n+\left(n-1\right)\sum_{i=1}^{t}c_{i}\right)-\left(C\cdot c_{0}+c_{0}^{2}+\sum_{0\le i<j\le t}c_{i}c_{j}\right)\left(n-1\right)\\
 & =C\left(C+S-\left(n-2\right)c_{0}\right)-\left(n-1\right)\left(c_{0}^{2}+c_{0}S+T\right)\\
 & =C^{2}-\left(n-2\right)c_{0}C+\left(n-2\right)S^{2}-\left(n-1\right)\left(c_{0}^{2}+T\right).
\end{align*}
By Lemma \ref{lem:if A subset of B}\eqref{enu:eigenvalues of alpha_a},
$\left\Vert \alpha_{B}\right\Vert \le\left|B\right|$ for every $B\subseteq\left[n\right]$,
so the right hand side of \eqref{eq:constants on the left} satisfies
that its $\ell_{2}$-norm is at most
\begin{align*}
\left\Vert \mathrm{RHS}\right\Vert  & \le\left(n-2\right)\sum_{i=1}^{t}c_{i}\left(C+c_{i}+c_{0}\right)+\left(n-3\right)\sum_{1\le i<j\le t}c_{i}c_{j}.\\
 & =\left(n-2\right)S\left(C+c_{0}\right)+\left(n-2\right)Q+\left(n-3\right)T\\
 & =\left(C-\left(n-1\right)c_{0}\right)\left(C+c_{0}\right)+\left(n-2\right)Q+\left(n-3\right)T\\
 & =C^{2}-\left(n-2\right)c_{0}C-\left(n-1\right)c_{0}^{2}+\left(n-2\right)Q+\left(n-3\right)T.
\end{align*}
So 
\begin{align*}
\mathrm{LHS}-\left\Vert \mathrm{RHS}\right\Vert  & \ge\left(n-2\right)S^{2}-\left(n-1\right)T-\left(n-2\right)Q-\left(n-3\right)T\\
 & =\left(n-2\right)\left[S^{2}-Q-2T\right]=0.\qedhere
\end{align*}
\end{proof}

\section{Further results around Caputo's conjecture for hypergraphs\label{sec:Further-results-around-Caputo-conjecture}}

Recall that Caputo's Conjecture \ref{conj:Caputo's--shuffles-conjecture}
concerns arbitrary finite weighted hypergraphs, and Theorem \ref{thm:triangles and the like}
proves only a special case of this conjecture. In the current section
we list a few further cases of this conjecture that we can prove. 

\subsection*{Tree-like hypergraphs}

Say that a hyper-edge in a hypergraph is a leaf if exactly one of
its vertices intersects other hyper-edges. We say that a finite hypergraph
is \emph{tree-like} if it is connected and after gradually removing
all leaves, one is left with a single hyper-edge\footnote{There are various other definitions for what the right higher-dimensional
analog of a tree is. Our definition here is by no means a conventional
one.}.

Given the results in this paper, it is easy to see that weighted tree-like
hypergraphs satisfy Caputo's conjecture. Indeed, use induction on
the number of vertices. The base case is trivial. Now in a general
hypergraph, pick any hyper-edge $E$ which is a leaf. If $E$ is a
singleton, removing it would only shift all spectra by a constant,
so assume that $\left|E\right|\ge2$ and without loss of generality
assume that $E=A\cup\left\{ 1\right\} $ and that the vertex $1$
is not contained in any other hyper-edge. 

As a special case of each of Theorems \ref{thm:squid - disjoint sets},
\ref{thm:squid - sets of size k containing a common (k-1)-subset}
or \ref{thm:squid - sets of cosize 1 in A0}, we get
\[
\left|A\right|\alpha_{A\cup\left\{ 1\right\} }\ge\left(\left|A\right|+1\right)\alpha_{A},
\]
so we can replace the hyper-edge $A\cup\left\{ 1\right\} $ with some
weight $w$, with the hyper-edge $A$ with weight $\frac{\left|A\right|+1}{\left|A\right|}w$,
and use the same argument as in \eqref{eq:CLR proof summary} and
\eqref{eq:proof-of-triangles} to use the induction hypothesis.

This argument can be used more generally whenever there is a hyper-edge
$E$ with some unique vertex $v$: a vertex belonging to no other
hyper-edge. We can always reduce to the case where $E$ is replaced
with $E\setminus\left\{ v\right\} $. (This simple trick was also
observed in \cite[\S4]{piras2010generalizations}.)

\subsection*{When all hyper-edges are of size $\ge n-1$}

The statement of Conjecture \ref{conj:Caputo's--shuffles-conjecture-Rep-version}
is equivalent to that for every non-negative weights $\left\{ w_{A}\right\} _{A\subseteq\left[n\right]}$,
the element $S_{G}=\sum_{A\subseteq\left[n\right]}w_{A}\alpha_{A}$
satisfies $\lambda_{\min}\left(S_{G},V_{\mu}\right)\ge\lambda_{\min}\left(S_{G},V_{\left(n-1,1\right)}\right)$
for all $\mu\vdash n$ with $\mu_{1}\le n-2$ (see the beginning of
Section \ref{subsec:A-general-inequality}). Now assume that $w_{A}=0$
whenever $\left|A\right|\le n-2$. By Corollary \ref{cor:alpha is constant when many squares outside first row},
\[
S_{G}=_{V_{\mu}}\sum_{A\subseteq[n]}w_{A}\left|A\right|
\]
for all $\mu\vdash n$ with $\mu_{1}\le n-2$. On the other hand,
\[
\left\Vert S_{G}|_{V_{\left(n-1,1\right)}}\right\Vert \le\sum_{A\subseteq\left[n\right]}w_{A}\left\Vert \alpha_{A}\right\Vert \le\sum_{A\subseteq\left[n\right]}w_{A}\left|A\right|
\]
by Lemma \ref{lem:if A subset of B}. We note that although it seems
related, this result does not seem to follow from Theorem \ref{thm:squid - sets of cosize 1 in A0}.

\subsection*{Small values of $n$}

Finally, let us also mention that using a combination of the results
of this paper and ideas taken from \cite{aldousorder} we managed
to prove Caputo's Conjecture \ref{conj:Caputo's--shuffles-conjecture-Rep-version}
for every hypergraph on at most six vertices. We do not include this
argument in the paper as we were not able to generalize it to larger
values of $n$.

\section{Some counterexamples\label{sec:Counterexamples}}

We elaborate here two concrete examples illustrating the failure of
the general inequality \eqref{eq:squid} to hold in all cases. 
\begin{example}
The following example was brought to our attention by P.~Caputo.
Consider two subsets $A_{1}=$$\left\{ 2\right\} $ and $A_{2}=\left\{ 2,3,4\right\} $
of $\left[4\right]=\left\{ 1,2,3,4\right\} $ with weights $c_{1}=c_{2}=1$.
The inequality \eqref{eq:squid} says in this case that 
\begin{equation}
4\left(\alpha_{\left\{ 1,2\right\} }+\alpha_{\left[4\right]}\right)\stackrel{?}{\ge}5\alpha_{\left\{ 2,3,4\right\} }+\left(\alpha_{\left\{ 2,3,4\right\} }-\alpha_{\left\{ 3,4\right\} }\right),\label{eq:squid fails}
\end{equation}
namely that $U\defi4\alpha_{\left\{ 1,2\right\} }+4\alpha_{\left[4\right]}-6\alpha_{\left\{ 2,3,4\right\} }+\alpha_{\left\{ 3,4\right\} }\stackrel{?}{\ge}0$.
However, direct computations show that $U$ has a negative value of
$-2$ in the representation $V_{\left(2,2\right)}$ (it is non-negative
in all other irreducible representations of $\Sym_{4}$). 

This failure of this inequality with these sets also shows that the
entire strategy laid out in Section \ref{sec:The-proof-of-Aldous}
fails for certain hypergraphs. For example, it fails for the hypergraph
$G$ on four vertices with positive weights exactly to the hyper-edges
$\left\{ 1,2\right\} $, $\left\{ 3,4\right\} $ and $\left\{ 1,2,3,4\right\} $.
If these weights are precisely $1$ each, then the failure of \eqref{eq:squid fails}
shows we cannot use the induction scheme of Section \ref{sec:The-proof-of-Aldous}
on any of the four vertices of $G$. This counterexample is a barrier
for all kinds of Caputo-Liggett-Righthammer-like approaches to the
conjecture.
\end{example}

\begin{example}
\label{exa:squid fails even with positive coefs}The following example,
found using computer simulations, illustrates that the general inequality
\eqref{eq:squid} may fail even when every non-vanishing $\alpha_{A}$
on the right hand side of \eqref{eq:squid} has a positive coefficient
(and see Remark \ref{rem:positive coefficient on the RHS}). Let $n=6$
and consider three subsets $A_{1}=\left\{ 2,3,4,5,6\right\} $, $A_{2}=\left\{ 2,3\right\} $
and $A_{4}=\left\{ 4,5,6\right\} $, with weights $c_{1}=\frac{3}{20},$$c_{2}=\frac{1}{10}$
and $c_{3}=\frac{1}{60}$, respectively. In this case, \eqref{eq:squid}
becomes
\[
\sum_{i=1}^{3}c_{i}\alpha_{A_{i}\cup\left\{ 1\right\} }\stackrel{?}{\ge}\frac{23}{120}\cdot\alpha_{\left\{ 2,3,4,5,6\right\} }+\frac{29}{240}\cdot\alpha_{\left\{ 2,3\right\} }+\frac{1}{360}\cdot\alpha_{\left\{ 4,5,6\right\} }.
\]
But in the irreducible representation $V_{\left(4,2\right)}$ the
difference between the left hand side and the right hand side has
a negative eigenvalue of roughly $-0.00517$, so this inequality does
not hold.
\end{example}

\section*{Appendix}

\begin{appendices}

\section{Code completing the proofs from Section \ref{subsec:hard cases}
for almost every $k$\label{sec:code hard cases for almost all k}}

Appendices \ref{sec:code hard cases for almost all k} and \ref{sec:All-remaining-cases small k}
contain the sagemath code used to finish the proof of Theorem \ref{thm:squid - sets of size k containing a common (k-1)-subset},
as described in Section \ref{sec:sets with large intersection}. In
particular, Section \ref{subsec:hard cases} discusses the non-negativity
of $\Gamma_{ij\ell m}$ on $V_{\mu}$ when $\mu\vdash k+4$ and $\mu_{1}=k+1$
or $\mu_{1}=k+2$, as well as the non-negativity of $\Gamma_{iij\ell}$
and of $\Gamma_{0ij\ell}$ when $\mu\vdash k+3$ and $\mu_{1}=k+1$.
That Section \ref{subsec:hard cases} explains how we verify the case
of $\Gamma_{\ijlm}$ for $k\ge4$, and we include the corresponding
code in the current Appendix \ref{sec:code hard cases for almost all k}.
We also include in the current appendix the code verifying the non-negativity
of $\Gamma_{iij\ell}$ and of $\Gamma_{0ij\ell}$ in the above representations
for every $k\ge3$. The remaining cases --- $\Gamma_{\ijlm}$ when
$k=2,3$ and $\Gamma_{iij\ell}$ and $\Gamma_{0ij\ell}$ when $k=2$
-- are included in Appendix \ref{sec:All-remaining-cases small k}.

The following code is used in the next two subsections. \\

\begin{lstlisting}
var('t') # for the characteristic polynomial

# accumulate 'value' into the (row,column) entry in the matrix mtrx
def accumulate( mtrx, row, column, value ):
	mtrx[basis.index(row),basis.index(column)] += value

# S is a subset of {1,...,5} (representing {1,k+1,k+2,k+3,k+4}
# b1 and b2 are two of the basis elements 
# This routine checks if the elements in the support of b1 are in the same orbit as 
# those in the support of b2, via the action of Sym(A_0 \cup S).
def in_orbit( S, b1, b2 ):
	c1 = list(b1)
	c2 = list(b2)
	tuple_len = len(b1) # the length of the tuple in every basis element (3 for ijlm, 2 for iijl and ijl0)
	for i in range(tuple_len):
		if c1[i] in S:
			c1[i] = 0  # the elements of (A_0 \cup S) are all in the same orbit and considered as "a star"
		if c2[i] in S:
			c2[i] = 0
	return c1==c2

# S is a subset of {1,...,5}, representing {1,k+1,k+2,k+3,k+4} (in practice, S will
# be of size one or two).
# B = A_0 \cup S  (so B is alpha_{A_i} or alpha_{A_ij})
# Construct the matrix giving the action of alpha_B on our space w.r.t. our basis, 
# analogously to Lemma 5.10
def alpha_B_matrix(S):
	rc = FracMatrix() # creates a zero matrix (136*136 or 21*21)
	size_B = len(S) + k - 1
	# As in (5.11), first add |B|*id
	for b in basis:
		accumulate(rc, b, b, size_B)
	# Now add the remaining part of (5.11)
	for b in basis:
		n_stars = b.count(0)
		n_B     = b.count(0) + sum([b.count(x) for x in S])
		coef = size_B
		for x in range(n_stars):
			coef *= (k - 1 - x)
		for x in range(n_B):
			coef /= (size_B - x)
		for b2 in basis:
			if in_orbit(S,b,b2):
				accumulate( rc, b2, b, -coef )

	return rc
  
# check the previous function: we take the weighted sum of every column of alpha_B
# to make sure we get zero
def test_alpha_B(S):
	V = MatrixSpace(C,1,dim_of_mtrx,sparse=False)
	v = V()
	for i in range(len(basis)):
		b = basis[i]
		n_stars = b.count(0)
		v[0,i] = falling_factorial(k-1,n_stars)
	return v*alpha_B_matrix(S)

# The Jordan product of two elements
def Jordan(x,y):
    return (x*y + y*x)/2
	
# return the value of w_ef as in (5.5)
# Assume e,f in {0,1,2,..}
def w(e,f):
	e,f = min(e,f),max(e,f)
	if e==0 and f==0: return (k-1)*k/2
	if 0<e and e==f: return k*(k+1)/2
	if e==0 and e<f: return k**2-1
	if 0<e and e<f: return k**2+k-1

# return the matrix representing the action of U_ef on our space, w.r.t. our basis. 
# We rely on (5.4)
# e,f assumed to belong to {0,2,3,4,5}
# Note that alpha_A_0 vanishes as we work in the kernel of alpha_A_0
def U(e,f):
	e,f = min(e,f),max(e,f)
	if e==0 and f==0: return (k-1)*alpha_B_matrix([1])
	if 0<e and e==f: return k*alpha_B_matrix([1,e]) - (k+1)*alpha_B_matrix([e])
	if e==0 and e<f: return k*alpha_B_matrix([1]) + (k-1)*alpha_B_matrix([1,f]) - k*alpha_B_matrix([f])
	if 0<e and e<f: return k*( alpha_B_matrix([1,e]) + alpha_B_matrix([1,f]) - alpha_B_matrix([e]) - alpha_B_matrix([f]) ) - alpha_B_matrix([e,f])

# Compute Lambda as in (5.8)
def Lambda(e1,f1,e2,f2):
    return FracMatrix() + w(e1,f1)*w(e2,f2) - Jordan( U(e1,f1)-w(e1,f1), U(e2,f2)-w(e2,f2) )




\end{lstlisting}

\subsection{The non-negativity of $\Gamma_{ij\ell m}$ on $W_{\left(k+1,1,1,1\right)}$
for $k\ge4$\label{subsec:code for Gamma ijlm}}

The way we verify that $\Gamma_{\ijlm}$ is non-negative on $W_{\left(k+1,1,1,1\right)}$
for every $k\ge4$ is detailed in Section \ref{subsec:hard cases}.
As described there, we work in $\Sym_{k+4}$ with $A_{0}=\left\{ 2,\ldots,k\right\} $,
$A_{i}=A_{0}\cup\left\{ k+1\right\} $ and so on. The 136 elements
of the basis we construct for $W_{\left(k+1,1,1,1\right)}\cap\ker\alpha_{A_{0}}$
are represented by triples of symbols from $\left\{ 1,k+1,k+2,k+3,k+4,*\right\} $.
In the code below we use $0$ instead of a star (``{*}''), and $2,\ldots,5$
instead of $k+1,\ldots,k+4$, respectively, (The number $1$ still
represents $1$.) The code below also uses the functions from the
beginning of this Section \ref{sec:code hard cases for almost all k}.
We remark that this piece of code took about two hours to complete
in SageMath 9.2 on a Linux machine.\\

\begin{lstlisting}
basis = [(x,y,z) for x in range(6) for y in range(6) for z in range(6) if x!=y and x!=z and y!=z] # no stars or one star
for x in range(1,6):
    basis += [(x,0,0),(0,x,0),(0,0,x)]  # two stars
basis += [(0,0,0)]          # three stars

B = PolynomialRing(QQ,['k'])  
C = FractionField(B)
C.inject_variables()
dim_of_mtrx = len(basis)
PolyMatrix = MatrixSpace(B,dim_of_mtrx,dim_of_mtrx,sparse=False)  # 136*136 matrices over Q[k]
FracMatrix = MatrixSpace(C,dim_of_mtrx,dim_of_mtrx,sparse=False)  # 136*136 matrices over Q(k)

# Compute Gamma_ijlm  as in (5.7)
def Gamma_matrix():
    rc = FracMatrix()
    rc += Lambda(2,3,4,5)
    rc += Lambda(2,4,3,5)
    rc += Lambda(2,5,3,4)
    rc += Lambda(3,4,2,5)
    rc += Lambda(3,5,2,4)
    rc += Lambda(4,5,2,3)
    return rc
   
#  This code is discussed at the very end of the proof of Theorem 1.9 in 
# Section 5.3: we verify that Gamma_ijlm is non-negative in W_(k+1,1,1,1) when k>=4
def verify_Thm_1_9_on_ijlm():
	G = Gamma_matrix()
	print( "G created")
	# By Lemma 5.11, the entries of (k^2 * G) are all polynomial in k
	H = PolyMatrix(k**2 * G)
	print("H created")
	p = H.characteristic_polynomial('t')
	print("p created")
	# check positive roots for k>=4
	q = p(k=k+3)(t=-t) * (-1)**p.degree(t)
	for i in range(q.degree(t)+1):
		for j in range(q.degree(k)+1):
			if q.coefficient(t**i * k**j) < 0:
				print( i,j,"ERROR" )
				return False                    
	print("All coefficients are non-negative, meaning that Gamma_ijlm is non-positive for every k>=4.")
	return True
\end{lstlisting}

\subsection{The non-negativity of $\Gamma_{iij\ell}$ and of $\Gamma_{0ij\ell}$
on $W_{\left(k+1,1,1\right)}$ for $k\ge3$\label{subsec:code for Gamma iijl and 0ijl}}

We now prove that $\Gamma_{iij\ell}$ and $\Gamma_{0ij\ell}$ ($i,j$
and $\ell$ are distinct positive integers) are both non-negative
on the reducible representation $W_{\left(k+1,1,1\right)}$ as explain
in Section \ref{subsec:hard cases}. We work in $\Sym_{k+3}$ and
without loss of generality take $A_{0}=\left\{ 2,\ldots,k\right\} $,
$A_{i}=A_{0}\cup\left\{ k+1\right\} $, $A_{j}=A_{0}\cup\left\{ k+2\right\} $
and $A_{\ell}=A_{0}\cup\left\{ k+3\right\} $. Similarly to the $\Gamma_{\ijlm}$
case, we may restrict our attention to 
\[
V_{k}\defi W_{\left(k+1,1,1\right)}\cap\ker\alpha_{A_{0}}.
\]
The elements of this subspace are invariant under every $\sigma\in\Sym_{A_{0}}$.
Our basis for $V_{k}$ consists of the following 21 elements:
\begin{itemize}
\item $v_{x,y}\defi\left(x,y\right)$ for distinct $x,y\in\left\{ 1,k+1,k+2,k+3\right\} $ 
\item ${\displaystyle v_{x,*}\defi\sum_{s\in\left\{ 2,\ldots,k\right\} }\left(x,s\right)}$
and similarly $v_{*,x}$ for $x\in\left\{ 1,k+1,k+2,k+3\right\} $
\item ${\displaystyle v_{*,*}\defi\sum_{\substack{s,t\in\left\{ 2,\ldots,k\right\} \colon\\
s\ne t
}
}\left(s,t\right)}$
\end{itemize}
Note that $v_{*,*}$ vanishes for $k=2$, so the current analysis
applies to $k\ge3$ only. (Recall that we assume throughout Section
\ref{sec:sets with large intersection} that $k\ge2$.) In the code
below we use $0$ instead of a star (\textquotedbl{*}\textquotedbl )
and $2,3,4$ instead of $k+1,k+2,k+3$, respectively. (The number
$1$ still represents $1$.) The code below also uses the functions
from the beginning of this Section \ref{sec:code hard cases for almost all k}.\\

\begin{lstlisting}
basis = [(i,j) for i in range(5) for j in range(5) if i!=j] # no stars or one star
basis.append((0,0)) # two stars

B = PolynomialRing(QQ,['k'])
C = FractionField(B)
C.inject_variables()
dim_of_mtrx = len(basis)
PolyMatrix = MatrixSpace(B,dim_of_mtrx,dim_of_mtrx,sparse=False)  # 21*21 matrices over Q[k]
FracMatrix = MatrixSpace(C,dim_of_mtrx,dim_of_mtrx,sparse=False)  # 21*21 matrices over Q(k)

# Compute Gamma_iijl  as in (5.7) 
def Gamma_iijl_matrix():
    rc = FracMatrix()
    rc += Lambda(2,2,3,4)
    rc += Lambda(2,3,2,4)
    rc += Lambda(2,4,2,3)
    rc += Lambda(3,4,2,2)
    return rc

# Compute Gamma_ijl0  as in (5.7) 
def Gamma_ijl0_matrix():
    rc = FracMatrix()
    rc += Lambda(0,2,3,4)
    rc += Lambda(0,3,2,4) 
    rc += Lambda(0,4,2,3)
    rc += Lambda(2,3,0,4)
    rc += Lambda(2,4,0,3)
    rc += Lambda(3,4,0,2)
    return rc

# The following scripts shows that the characteristic polynomial of Gamma matrix
# has no negative roots for k>=3. We construct the characteristic polynomial p of 
# the (normalized) matrix associated with the action of Gamma on our space.  
# This is a polynomial in the variable t with coefficients that are polynomials 
# in k over the integers. We verify that for every value k>=3 this polynomial has 
# only non-negative roots. 
# We use the easy criterion that a real-rooted polynomial with real coefficients 
# has only non-negative roots if and only if the signs of coefficients alternate 
# between non-negative and non-positive. By substituting t=-t and multiplying by
# (-1)^degree_of_poly we want to show that the coefficients are positive for all 
# k>=3. By substituting k=k+2, we get p', and it is enough to verify that the 
# coefficients of p' are now are polynomials in k with non-negative coefficients. 
def verify_Thm_1_9_on_iijl_or_ijl0():
	for case in ['iijl','0ijl']:
		if case=="iijl": G = Gamma_iijl_matrix()
		if case=="0ijl": G = Gamma_ijl0_matrix()
		# By Lemma 5.11, the entries of (k^2 * G) are all polynomial in k
		H = PolyMatrix(k**2 * G)
		p = H.characteristic_polynomial('t')
		q = p(k=k+2)(t=-t) * (-1)**p.degree(t)
		for i in range(q.degree(t)+1):
			for j in range(q.degree(k)+1):
				if q.coefficient(t**i * k**j) < 0:
					print( i,j,"ERROR" )
					return False
		print( "In case " + case + " all coefficients are non-negative so Gamma_" + case + " is non-negative for k>=3" )
	return True
\end{lstlisting}

\section{Code for all remaining, small $k$ cases from Section \ref{sec:sets with large intersection}
\label{sec:All-remaining-cases small k}}

Finally, the following code deals with certain remaining cases from
the proof of Theorem \ref{thm:squid - sets of size k containing a common (k-1)-subset}
in Section \ref{sec:sets with large intersection}. It includes a
function to check the cases $k=2,3,4,5$ in Lemma \ref{lem:mu1=00003Dk 00ij},
$k=2$ in Lemma \ref{lem:mu1=00003Dk 0iij } and $k=2$ in Lemma \ref{lem:mu1=00003Dk 0ijl}.
It also verifies the non-negativity of $\Gamma_{\ijlm}$ on the irreducible
representations discussed in Section \ref{subsec:hard cases} when
$k=2,3$, and the non-negativity of $\Gamma_{iij\ell}$ and $\Gamma_{0ij\ell}$
in the irreducible representations discussed in Section \ref{subsec:hard cases}
when $k=2$.\\

\begin{lstlisting}
# for a given subset A of [n], return alpha_A -- the sum of (1-p) for all
# permutations p supported on A, divided by (|A|-1)!. 
# subset is assumed to be a list
def alpha_of_set( subset, n ):
    Sn = SymmetricGroup(n)
    Sk = SymmetricGroup(subset)
    B = GroupAlgebra(Sn,QQ)
    m = len(subset)
    if m<=1:
        return B.zero()
    perms = [Sn(p) for p in Sk]
    rc = B.zero()
    for p in perms:
        rc += B(Sn.one()) - B(p)
    rc = rc/factorial(m-1)
    return rc

# elem - element in the group ring R[S_n]
# p - a partition of n (a list)
# test if the image of the element through the irreducible representation of S_n
# corresponding to p is non-negative.
# We use the easy criterion that a real-rooted polynomial with real coefficients
# has only non-negative roots if and only if the signs of coefficients alternate
# between non-negative and non-positive.
def is_nonnegative( elem, p, verbose=False ):
    Sn = elem.parent().group()
    spc = SymmetricGroupRepresentation(p) #'specht' realization
    A = sum([spc.representation_matrix(Permutation(g))*elem.coefficient(g) for g in Sn if elem.coefficient(g)!=0])
    pol = A.characteristic_polynomial()
    if verbose: print ("    The characteristic polynomial p(x) is:", pol)
    pol2 = pol(x=-x) * (-1)**pol.degree()
    if verbose: print ("                     up to a sign  p(-x) =", pol2)
    for i in range(pol.degree()):
        if pol2.coefficient(x,i) < 0:
            if verbose: print("      There are negative roots. This element is NOT non-negative.")
            return False
    if verbose: print ("    all the coefficients are nonnegative, so all the roots of p(x) are nonnegative.")
    return true

# The Jordan product of two elements
def Jordan(x,y):
    return (x*y + y*x)/2

# A_0 = {2,3,...,k}
def A0(k):
	return list(range(2,k+1))
# A_0 \cup {1} = {1,2,3,...,k}
def A01(k):
	return list(range(1,k+1))
# A_i = {2,3,...,k} \cup {k+i}   we assume i=1,2,3,... 
def Ai(k,i):
	rc = A0(k)
	rc.append(k+i)
	return rc
# A_i \cup {1} = {1,2,3,...,k} \cup {k+i}
def Ai1(k,i):
	rc = Ai(k,i)
	rc.append(1)
	return rc
# A_i \cup A_j = {2,3,...,k} \cup {k+i,k+j}
def Aij(k,i,j):
	if i==j:
		print ("ERROR")
		return Ai(k,i)
	rc = Ai(k,i)
	rc.append(k+j)
	return rc

# U_ij inside S_n, as defined in Equation (5.4)
# We assume n is large enough (at least k+max(i,j))
def U(k,i,j,n):
	i,j = min(i,j),max(i,j)  # sort i,j
	if i==0 and j==0: return (k-1)*alpha_of_set(A01(k),n) - k*alpha_of_set(A0(k),n)
	if i>0  and i==j: return k*alpha_of_set(Ai1(k,i),n) - (k+1)*alpha_of_set(Ai(k,i),n)
	if i==0 and i<j:  return k*alpha_of_set(A01(k),n) + (k-1)*alpha_of_set(Ai1(k,j),n) - (k+1)*alpha_of_set(A0(k),n) - k*alpha_of_set(Ai(k,j),n)
	if i>0  and i<j:  
		return k*(alpha_of_set(Ai1(k,i),n) + alpha_of_set(Ai1(k,j),n) - alpha_of_set(Ai(k,i),n) - alpha_of_set(Ai(k,j),n)) - alpha_of_set(Aij(k,i,j),n) - alpha_of_set(A0(k),n)

# the numbers w_ij, as defined in Equation (5.5)
def w(k,i,j):
	i,j = min(i,j),max(i,j)  # sort i,j
	if i==0 and j==0: return ZZ(k*(k-1)/2)
	if i>0  and i==j: return ZZ(k*(k+1)/2)
	if i==0 and i<j:  return ZZ(k*k - 1)
	if i>0  and i<j:  return ZZ(k*k + k - 1)

# Lambda as in Equation (5.8) in the group ring of S_n
def Lambda(k,e1,f1,e2,f2,n):
	return w(k,e1,f1)*w(k,e2,f2) - Jordan( U(k,e1,f1,n)-w(k,e1,f1), U(k,e2,f2,n)-w(k,e2,f2) )

# test specific "Gamma's"
def Gamma_00ij(k):
	n = k+2
	return Lambda(k,0,0,1,2,n) + Lambda(k,1,2,0,0,n) + Lambda(k,0,1,0,2,n) + Lambda(k,0,2,0,1,n)
def Gamma_0iij(k):
	n = k+2
	return Lambda(k,1,1,0,2,n) + Lambda(k,0,2,1,1,n) + Lambda(k,0,1,1,2,n) + Lambda(k,1,2,0,1,n)    	
def Gamma_0ijl(k):
	n = k+3
	return Lambda(k,0,1,2,3,n) + Lambda(k,0,2,1,3,n) + Lambda(k,0,3,1,2,n) + Lambda(k,1,2,0,3,n) + Lambda(k,1,3,0,2,n) + Lambda(k,2,3,0,1,n)   
def Gamma_iijl(k):
	n = k+3
	return Lambda(k,1,1,2,3,n) + Lambda(k,1,2,1,3,n) + Lambda(k,1,3,1,2,n) + Lambda(k,2,3,1,1,n)
def Gamma_ijlm(k):
	n = k+4
	return Lambda(k,1,2,3,4,n) + Lambda(k,1,3,2,4,n) + Lambda(k,1,4,2,3,n) + Lambda(k,2,3,1,4,n) + Lambda(k,2,4,1,3,n) + Lambda(k,3,4,1,2,n)

# Verify that Lemma 5.7 holds for k=2,3,4,5
def verify_Lemma_5_7():
    for k in range(2,6):
        print("Testing the lemma for k =",k)
        g = Gamma_00ij(k)
        for p in [[k,2],[k,1,1]]: # These are the partitions mu of (k+2)  with mu_1=k
            print("  Testing the representation", p)
            if is_nonnegative( g, p, verbose=True ): print("    so far so good")
            else: print("Oh no! The paper is doomed"); return
    print("Hooray! Lemma 5.7 is verified")
    
# Verify that Lemma 5.8 holds for k=2
def verify_Lemma_5_8():
    k=2
    print("Testing the lemma for k=",k)
    g = Gamma_0iij(k)
    for p in [[k,2],[k,1,1]]:
        print("  Testing the representation", p)
        if is_nonnegative( g, p, verbose=True ): print("    so far so good")
        else: print("Oh no! The paper is doomed"); return
    print("Hooray! Lemma 5.8 is verified")

# Verify that Lemma 5.9 holds for k=2
def verify_Lemma_5_9():
    k=2
    print("Testing the lemma for k=",k)
    g = Gamma_0ijl(k)
    for p in [[k,2,1],[k,1,1,1]]: # When n=k+3=5, there are only two partitions with first row of length 2
        print("  Testing the representation", p)
        if is_nonnegative( g, p, verbose=True ): print("     so far so good")
        else: print("Oh no! The paper is doomed"); return
    print("Hooray! Lemma 5.9 is verified")

# Verify that cases discussed in Section 5.3 hold for:
#  Gamma_ijlm and k=2,3
#  Gamma_iijl and k=2
#  Gamma_0ijl and k=2
def verify_remainig_cases_Section_5_3():
	# Test Gamma_ijlm
	for k in [2,3]:
		print("Testing Gamma_ijlm, k+1 <= mu_1 <= k+2, with k=",k)
		g = Gamma_ijlm(k)
		for p in [ [k+1,3], [k+1,2,1], [k+1,1,1,1], [k+2,2], [k+2,1,1] ]: # These are the partitions mu of k+4 with mu_1=k+1 or k+2
			print("  Testing the representation", p)
			if is_nonnegative( g, p, verbose=True ): print("    so far so good")
			else: print("Oh no! The paper is doomed"); return
	# Now test Gamma_iijl 
	k=2
	print("Testing Gamma_iijl, with k=",k," and mu_1=k+1")
	g = Gamma_iijl(k)
	for p in [ [k+1,2], [k+1,1,1] ]: # These are the partitions mu of k+3 with mu_1=k+1
		print("  Testing the representation", p)
		if is_nonnegative( g, p, verbose=True ): print("    so far so good")
		else: print("Oh no! The paper is doomed"); return
	# Now test Gamma_0ijl 
	k=2
	print("Testing Gamma_0ijl, with k=",k," and mu_1=k+1")
	g = Gamma_0ijl(k)
	for p in [ [k+1,2], [k+1,1,1] ]: # These are the partitions mu of k+3 with mu_1=k+1
		print("  Testing the representation", p)
		if is_nonnegative( g, p, verbose=True ): print("    so far so good")
		else: print("Oh no! The paper is doomed"); return
	print("Hooray! All sporadic cases discussed in Section 5.3 are verified")
\end{lstlisting}

\end{appendices}

\bibliographystyle{alpha}
\bibliography{GenOctopus-bib}

\newcommand{\etalchar}[1]{$^{#1}$}
\begin{thebibliography}{ACD{\etalchar{+}}20}

\bibitem[ACD{\etalchar{+}}20]{aldous2020life}
D.~Aldous, P.~Caputo, R.~Durrett, A.~E. Holroyd, P.u Jung, and A.~L. Puha.
\newblock The life and mathematical legacy of {T}homas {M}. {L}iggett.
\newblock {\em Notices Amer. Math. Soc.}, 68(1), 2020.

\bibitem[AK13]{aldousorder}
G.~Alon and G.~Kozma.
\newblock Ordering the representations of {$S_n$} using the interchange
  process.
\newblock {\em Canad. Math. Bull.}, 56(1):13--30, 2013.

\bibitem[AK20]{alon2020comparing}
G.~Alon and G.~Kozma.
\newblock Comparing with octopi.
\newblock {\em Ann. Inst. H. Poincar\'{e} Probab. Statist.}, 56(4):2672--2685,
  2020.

\bibitem[BC24]{bristiel2024entropy}
A.~Bristiel and P.~Caputo.
\newblock Entropy inequalities for random walks and permutations.
\newblock {\em Ann. Ins. Henri Poincar{\'e} (B) Probab. Stat.}, 60(1):54--81,
  2024.

\bibitem[Ces16]{cesi2016few}
F.~Cesi.
\newblock A few remarks on the octopus inequality and {A}ldous’ spectral gap
  conjecture.
\newblock {\em Comm. Algebra}, 44(1):279--302, 2016.

\bibitem[Che17]{chen2017moving}
J.~P. Chen.
\newblock The moving particle lemma for the exclusion process on a weighted
  graph.
\newblock {\em Electron. Commun. Probab.}, 22:1--13, 2017.

\bibitem[CLR10]{caputo2010proof}
P.~Caputo, T.~Liggett, and T.~Richthammer.
\newblock Proof of {A}ldous' spectral gap conjecture.
\newblock {\em J. Amer. Math. Soc.}, 23(3):831--851, 2010.

\bibitem[DS81]{diaconis1981generating}
P.~Diaconis and M.~Shahshahani.
\newblock Generating a random permutation with random transpositions.
\newblock {\em Z. Wahrscheinlichkeitstheorie verw. Gebiete}, 57(2):159--179,
  1981.

\bibitem[FH13]{fulton2013representation}
W.~Fulton and J.~Harris.
\newblock {\em Representation theory: a first course}, volume 129.
\newblock Springer Science \& Business Media, 2013.

\bibitem[HJ96]{handjani1996rate}
S.~Handjani and D.~Jungreis.
\newblock Rate of convergence for shuffling cards by transpositions.
\newblock {\em J. Theoret. Probab.}, 9:983--993, 1996.

\bibitem[MSS14]{marcus2014ramanujan}
A.~W. Marcus, D.~A. Spielman, and N.~Srivastava.
\newblock Ramanujan graphs and the solution of the {K}adison-{S}inger problem.
\newblock In {\em 2014 International Congress of Mathematicans, ICM 2014},
  volume~3, pages 363--386, 2014.

\bibitem[Pir10]{piras2010generalizations}
D.~Piras.
\newblock {\em Generalizations of {A}ldous’ Spectral Gap Conjecture}.
\newblock PhD thesis, Tesi di Laurea, Universit{\'a} degli Studi Roma Tre,
  2010.

\bibitem[PP20]{parzanchevski2020aldous}
O.~Parzanchevski and D.~Puder.
\newblock Aldous’s spectral gap conjecture for normal sets.
\newblock {\em Trans. Amer. Math. Soc.}, 373(10):7067--7086, 2020.

\end{thebibliography}

\noindent Gil Alon, Department of Mathematics and Computer Science,
The Open University of Israel, Raanana, 43107, Israel\\
\texttt{gilal@openu.ac.il}~\\
\texttt{}~\\
Gady Kozma, Faculty of Mathematics and Computer Science, Weizmann
Institute, Rehovot, 7610001, Israel\\
\texttt{gady.kozma@weizmann.ac.il}~\\

\noindent Doron Puder, School of Mathematical Sciences, Tel Aviv University,
Tel Aviv, 6997801, Israel\\
\texttt{doronpuder@gmail.com}
\end{document}